\documentclass[12pt]{article} 
\usepackage[margin=1in]{geometry}
\usepackage{url}
\usepackage[utf8]{inputenc}
\usepackage[english]{babel}
\usepackage{amsthm}
\usepackage{upgreek}
\usepackage{bbm, dsfont}

\usepackage[mathscr]{euscript} 
 \let\mathscr\relax 
\usepackage{cite,mathtools} 
\usepackage{latexsym}
\usepackage{amssymb,amsbsy,amsmath,amsfonts,amssymb,amscd,amsfonts,mathrsfs}
\usepackage{graphicx}
\usepackage{color} 

\usepackage{txfonts}
\usepackage{mathtools}

\DeclareMathOperator*{\argmax}{arg\,max}

\newcommand {\e} {\varepsilon}
\newtheorem{remark}{Remark}[]

\numberwithin{equation}{section}

\usepackage{titling}
\predate{\begin{flushleft}}
\postdate{\end{flushleft}}
\allowdisplaybreaks
\usepackage[font=rm, labelfont=bf]{caption}

\title{Derivation and travelling wave analysis of phenotype-structured haptotaxis models of cancer invasion} 

\author{\normalsize
Tommaso Lorenzi$^{1}$,
Fiona R Macfarlane$^{2,*}$,
Kevin J Painter$^{3}$
}

\date{$^{1}$ Department of Mathematical Sciences ``G. L. Lagrange'', Politecnico di Torino, 10129 Torino, Italy;\\
 $^{2}$ School of Mathematics and Statistics, University of St Andrews, United Kingdom;\\
 $^{3}$ Inter-university Department of Regional and Urban Studies and Planning, Politecnico di Torino, 10129 Torino, Italy.\\
 \ \\
 $^*$Corresponding author: frm3@st-andrews.ac.uk}

\begin{document}
\maketitle
\begin{abstract}
We formulate haptotaxis models of cancer invasion wherein the infiltrating cancer cells can occupy a spectrum of states in phenotype space, ranging from `fully mesenchymal’ to `fully epithelial’. The more mesenchymal cells are those that display stronger haptotaxis responses and have greater capacity to modify the extracellular matrix (ECM) through enhanced secretion of matrix-degrading enzymes (MDEs). However, as a trade-off, they have lower proliferative capacity than the more epithelial cells. The framework is multiscale in that we start with an individual-based model that tracks the dynamics of single cells, which is based on a branching random walk over a lattice representing both physical and phenotype space. We formally derive the corresponding continuum model, which takes the form of a coupled system comprising a partial integro-differential equation for the local cell population density function, a partial differential equation for the MDE concentration, and an infinite-dimensional ordinary differential equation for the ECM density. Despite the intricacy of the model, we show, through formal asymptotic techniques, that for certain parameter regimes it is possible to carry out a detailed travelling wave analysis and obtain invading fronts with spatial structuring of phenotypes. Precisely, the most mesenchymal cells dominate the leading edge of the invasion wave and the most epithelial (and most proliferative) dominate the rear, representing a bulk tumour population. As such, the model recapitulates similar observations into a front to back structuring of invasion waves into leader-type and follower-type cells, witnessed in an increasing number of experimental studies over recent years.
\end{abstract}

\section{Introduction}

\subsection{Biological background}
Phenomena of collective cell migration have received a significant amount of interest in recent years, and a particularly large literature has been devoted to their role during cancer invasion processes~\cite{friedl2012classifying}. Histopathological analysis of tissue specimens reveals a plurality of patterns within invading cell fronts, ranging from individually migrating cells to collective strands and clusters that infiltrate the surrounding healthy tissue. Increased invasiveness forms one of the key traits of cancer metastasis, and consequently presents a significant impediment to successful treatment. Understanding the underlying biological processes is therefore of manifest interest.

Invading cell fronts frequently present significant phenotypic heterogeneity, and particular attention has focussed on the extent to which a separation into `leader' and `follower' cells contributes to invasive spread~\cite{vilchez2021decoding}. Leader cells are those that (seemingly) drive the invasion process: positioned at or near the front, modifying the extracellular matrix (ECM) -- i.e. the network of macromolecules in which cells reside, and coordinating with follower cells to facilitate their collective invasion into healthy tissue. These leader cells can be derived from various sources, both from tumour cells that have undergone a phenotypic transition (e.g. following genetic mutations or epimutations) or from surrounding stromal cells, such as fibroblasts, that have been activated and co-opted through factors in the tumour microenvironment~\cite{vilchez2021decoding}.

For carcinomas -- cancers of epithelial origin -- tumour-derived leader cells will typically have undergone an epithelial to mesenchymal transition (EMT), i.e. undergone a loss of their epithelial nature and acquired mesenchymal characteristics~\cite{brabletz2018emt}. Epithelial cells are often tightly bound through strong cell-to-cell adhesion, therefore downregulation of adhesion can release leaders from the bonds that control their normal position in the tissue. Mesenchymal characteristics involve (significantly) enhanced motility and strong interactions with the ECM. These interactions include a capacity to modify the matrix and microenvironment, both mechanically and chemically, in a manner that can facilitate the invasion of other leader and follower cells. First, increased production of fibronectin can increase matrix adhesivity, allowing cells to gain more traction~\cite{konen2017image,schwager2019cell,vilchez2021decoding}; the directed movement that results from migration up a matrix adhesivity gradient is referred to as haptotaxis~\cite{carter1965principles}. Cells may also realign fibres, leading to an oriented matrix that directs invasion along certain paths~\cite{ray2021aligned}. Leader cells may also start to secrete (or increase the secretion of) matrix-degrading enzymes (MDEs)~\cite{chen2020makes,vilchezref78,schwager2019cell,vilchez2021decoding,vilchezref13}. In turn, this can reduce the volume fraction occupied by the ECM and liberate space into which cells can migrate (or proliferate). Beyond these physical alterations to the microenvironment, leader cells may also drive invasion through chemical means, e.g. through secreted factors leading to  chemoattractant gradients that direct follower cell movement~\cite{chen2020makes,vilchezref14,konen2017image,vilchez2021decoding}.

While leader cells may have enhanced motility and greater capacity to alter the surrounding microenvironment, trade-offs may arise in the form of reduced proliferative potential. Analyses into heterogeneous groups composed from distinct follower and leader subpopulations indicate that the former can be significantly more proliferative~\cite{konen2017image}. Furthermore, followers can produce factors that promote a level of growth within the leader cells, hence maintaining the balance between leaders and followers across the overall cell population~\cite{konen2017image}. Thus, collective invasion in certain leader-follower cancer cell populations may be a cooperative process, with each subpopulation playing a key role in promoting tumour expansion. 

Dichotomising invading cells into leader and follower subtypes can be notionally convenient, but masks the possibility that cells may fall into (significantly) more than two fundamental phenotypic states. Unlike the `regulated' EMTs that occur during embryogenesis or wound healing, EMT within cancer cells can be highly variable, ranging from `partial' to `full'~\cite{brabletz2018emt,chen2020makes,vilchezref27,vilchezref26,vilchezref24,vilchez2021decoding}. Partial here indicates a cell has only acquired a subset of mesenchymal characteristics, for example resulting in only a slight increase in motility, or a part reduction of proliferative potential. Moreover, it has been shown that cells can undergo changes to these phenotypes over time~\cite{vilchezref31,vilchez2021decoding,vilchezref30,vilchezref32}. Consequently, a more accurate picture of leader-follower heterogeneity in invading cancers would be that each cell occupies a fluid position within a quasi-continuous phenotype space. Investigating the role of trade-offs in leader-follower collective migration may form a beneficial treatment target. In fact, it has been shown that leader cells are generally more resistant to treatment~\cite{vilchezref15,vilchez2021decoding}; however, if leader cells are removed then the invasion of the tumour stops~\cite{konen2017image}. Furthermore, through inhibiting MDE secretion by leader cells then further invasion can be slowed or prevented~\cite{vilchezref23}.

\subsection{Mathematical modelling background}
The use of mathematical modelling to investigate the invasion processes of cancer cells forms a well developed area of research, and for more details we refer to the {extensive reviews detailed in~\cite{alfonso2017biology,franssen2019mathematical,sfakianakis2020mathematical}}. A significant part of the literature on this subject has focussed on haptotaxis-fuelled invasion, in which a cancer cell population secretes proteolytic factors that alter the surrounding matrix to generate adhesivity gradients. Early models in this field have been formulated as coupled systems of partial differential equations (PDEs) and infinite-dimensional ordinary differential equations (ODEs), which rely on the assumption that cell migration results from the superposition of random motion, modelled as linear diffusion, and haptotaxis, modelled as advection according to the ECM gradient, e.g.~\cite{anderson2000mathematical,perumpanani1999extracellular}. Focussing on a 1D spatial scenario, as a prototypical example of these haptotaxis models of cancer invasion we can consider the following system
\begin{equation}
	\begin{cases}
		\displaystyle{\partial_t \rho} - \partial_x \Big(D \, \partial_x \rho - \rho \, \mu \, \partial_x E\Big) = R(\rho) \, \rho,\\
		\ \\
		\displaystyle{\partial_t M} - D_M \, \partial^2_{xx} M = \gamma \, \rho - \kappa_M \, M,\\ \ \\
		\displaystyle{\partial_t E}=-\kappa_E \, E \, M,
	\end{cases}\label{eq:generic}
	\quad
	x \in \mathbb{R},
\end{equation}
where the functions $\rho \equiv \rho(t,x)$, $M \equiv M(t,x)$, and $E \equiv E(t,x)$ model, respectively, the cancer cell density, the MDE concentration, and the ECM density at time $t \in \mathbb{R}^+$ and position $x \in \mathbb{R}$. In the system~\eqref{eq:generic}, the parameters $D$ and $\mu$ are the random motility coefficient and the coefficient of haptotaxis sensitivity (i.e. sensitivity to matrix adhesivity gradients) of cancer cells, respectively, while the function $R(\rho)$ is the net growth rate of the cell population due to the proliferation (i.e. division and death) of cancer cells. The dependence of this function on cell density $\rho$ takes into account density-dependent inhibition of growth (i.e. the fact that cessation of cell division occurs once a critical value of the cell density is reached). Moreover, the parameter $D_M$ is the diffusion coefficient of MDEs, the parameter $\gamma$ is the rate of MDE production by cancer cells, and the parameter $\kappa_M$ is the rate of natural decay of MDEs. Finally, the parameter $\kappa_E$ is the rate at which the ECM is degraded by MDEs upon contact.

Haptotaxis models of cancer invasion of the form of~\eqref{eq:generic} have been studied analytically through proof of local existence, global existence, boundedness, and uniqueness of solutions~\cite{tao2010density}, investigations into blow-up of solutions~\cite{shangerganesh2019finite}, and analysis of 2D radially symmetric solutions~\cite{bortuli2021group}. Numerical simulations have also been used to investigate models of the form of \eqref{eq:generic}~\cite{ganesan2017galerkin}. Further models extend systems of this form, for example to include more detailed mechanisms of ECM remodelling and enzyme activities~\cite{andasari2011mathematical,chaplain2005mathematical,nguyen2018mathematical}, or through the inclusion of non-local terms to incorporate the effects of cell-cell adhesion~\cite{gerisch2008mathematical, painter2010impact}. Models have also been formulated to include detailed cell mechanics -- e.g. stress, strain, elasticity, adhesion, transport by velocity fields, and other interaction forces -- in the context of invasive melanoma growth through different skin layers~\cite{ciarletta2011radial}.

While model~\eqref{eq:generic} has been restricted to a homogeneous cancer population, cognate models have also been developed that explore the consequences of phenotypic heterogeneity on invasion. Binary state models consider a division of cancer cells into two phenotypic states. These include models in which cells of the invading population switch between a proliferating state and a migrating state, to investigate ramifications of the ``go-or-grow" hypothesis for glioma growth~\cite{kolbe2020modeling,pham2012density,stepien2018traveling}. Invasion models that feature two competing phenotypes with distinct migratory and proliferation properties have also been formulated in the context of acid-mediated invasion, where the heterogeneity extends to distinct acid-resistance and matrix-altering behaviour~\cite{strobl2020mix}. 

Greater phenotypic heterogeneity can be accounted for through the inclusion of more discrete states, but this becomes impractical if the phenotypic space becomes almost continuous. As such, an alternative approach is to extend a model like~\eqref{eq:generic} to include a continuous structuring variable representing intercellular variability in certain phenotypic characteristics~\cite{perthame2006transport}. In the context of cell invasion type dynamics, though not specifically in cancer, models of this nature have been developed to describe how a trade-off between chemotactic ability and proliferation may shape the phenotypic structuring of chemotaxis-driven growth processes~\cite{lorenzi2022trade}, and how trade-offs between mobility and proliferation may impact on density- or pressure-driven growth processes~\cite{lorenzi2021,macfarlane2022individual}; directly relevant to cancer, structured phenotype models of this type have also been developed to explore avascular tumour growth~\cite{fiandaca2022phenotype} and the evolutionary dynamics underpinning the emergence of intra-tumour phenotypic heterogeneity~\cite{fiandaca2021mathematical,lorenzi2018role,villa2021modeling}. {In recent work by Guilberteau \textit{et al.}~\cite{guilberteau2023integrative}, the authors presented a PIDE model that captures transitions between fully-epithelial, hybrid epithelial/mesenchymal, and fully-mesenchymal cell states, and demonstrated this model to be capable of reproducing experimental observations into the dynamics of EMT. This model, however, does not account for spatial dynamics or invasion processes of cells.}

While the above discussions have focussed on continuum models, which provide a population-level description of cell dynamics, it is important to note that a very large number of modelling studies have explored haptotaxis-driven cancer invasion using individual-based models (i.e. models that track the dynamics of single cells)~\cite{wang2015simulating,west2022agent}. Advantages lie in the ability of these models to capture the dynamics and stochasticity of single-cell movement, and a notable early example within the context of cancer invasion was developed in~\cite{anderson2000mathematical}. Here, a model of the type of~\eqref{eq:generic} was discretised in space, with the discretised terms subsequently used to specify probabilities of movement in different directions, according to the ECM distribution. Extensions were introduced in the model considered in~\cite{anderson2006tumor}, where each cell was allowed to undergo random movement, haptotaxis up the gradient of the ECM, produce MDEs, consume oxygen, and undergo cell cycle controlled proliferation depending on the availability of oxygen and free space. Moreover, this model comprised cells of different discrete phenotypic states, controlling aspects such as each cell's adhesion, oxygen consumption, haptotactic ability, secretion rate of MDEs, and proliferative potential. Recently, this underlying modelling framework has been further extended to investigate the role of two specific phenotypes, namely epithelial and mesenchymal phenotypes, in cancer invasion and metastasis~\cite{franssen2019mathematical}.

The original framework in the above method constituted of starting with a continuum model and subsequently discretising it to derive the governing rules for cell movement in a corresponding individual-based model~\cite{anderson2000mathematical}. An alternative approach for transitioning between discrete and continuum descriptions of cell motion is to first postulate a model at the single-cell level and then employ coarse-graining procedures to derive a continuous description; these methods have been extensively adopted in recent decades, particularly in the context of motivating PDE models to describe taxis-like behaviours, e.g.~\cite{painter2019,penington2011building,stevens1997aggregation}.

\subsection{Outline}
In this paper, we consider the following generalisation of the haptotaxis model of cancer invasion~\eqref{eq:generic}, where the continuous structuring variable $y \in [0,Y] \subset \mathbb{R}^+$, with $Y > 0$, represents the cell phenotype (i.e. the position of the cells in phenotypic space) and captures intercellular variability in haptotactic response, proliferative potential, and production of MDEs:
\begin{equation}
\label{eq:PDE}
\begin{cases}
\displaystyle{\partial_t n - \partial_x\Big(D \, \partial_x n - n \, \chi(y) \, \partial_x E\Big) = R(y,\rho) \, n + \lambda \, \partial_{yy} n},\\ \ \\
\displaystyle{\rho(t,x) := \int_0^Y n(t,x,y)\ \mathrm{d}y},\\ \ \\
\displaystyle{\partial_t M - D_M\partial_{xx}M = \int_0^Y p(y) \, n(t,x,y) \ \mathrm{d}y- \kappa_M M},\\ \ \\
\displaystyle{\partial_t E=-\kappa_E \, E \, M},
\end{cases}
\quad
(x,y) \in \mathbb{R} \times (0,Y).
\end{equation}
Compared to model~\eqref{eq:generic}, here the PDE for the cell density, $\rho(t,x)$, is replaced by the partial integro-differential equation (PIDE)~\eqref{eq:PDE}$_1$ for the local cell population density function, $n \equiv n(t,x,y)$, which is linked to the cell density through the relation~\eqref{eq:PDE}$_2$. Moreover, the functions $\chi(y)$, $p(y)$, and $R(y,\rho)$ are, respectively, the haptotaxis sensitivity coefficient, the MDE production rate, and the net growth rate of the cell population density due to the proliferation (i.e. division and death) of cancer cells with phenotype $y$. Finally, the diffusion term on the right-hand side of the PIDE~\eqref{eq:PDE}$_1$ models the effect of phenotypic changes, which occur at rate $\lambda$. 

We first formulate a phenotype-structured individual-based haptotaxis model of cancer invasion (cf. Section~\ref{sec:model}), where the dynamics of individual cells are governed by a set of rules that result in a branching random walk over a lattice~\cite{hughes1995random}, which represents both physical and phenotype spaces. In this model, the rules governing cell dynamics are coupled with a balance equation for the MDE concentration and a balance equation for the density of ECM. Then, using an extension of the limiting procedure that we previously employed in~\cite{bubba2020discrete,chaplain2020bridging,macfarlane2020hybrid,macfarlane2022individual}, we formally derive the model~\eqref{eq:PDE} as the continuum limit of this individual-based model (cf. Section~\ref{sec:contmodel} and Appendix~\ref{app:derivation}). After that, building upon the formal asymptotic method that we developed in~\cite{lorenzi2022trade,lorenzi2021}, we carry out travelling wave analysis of an appropriately rescaled version of the model~\eqref{eq:PDE} (cf. Section~\ref{sec:twanalysis}). The results obtained provide a mathematical formalisation for the idea that trade-offs between proliferative potential and the ability to sense spatial gradients of ECM and produce MDEs may promote the emergence of phenotypically structured invading cell fronts. Specifically, wherein leader cells (i.e. cells with a higher haptotactic and MDE production ability but a lower proliferative potential) are localised at the leading edge of the front, while follower cells (i.e. cells with a higher proliferative potential but a lower haptotactic and MDE production ability) occupy the region behind the leading edge. Finally, we report on numerical solutions of such a rescaled continuum model and numerical simulations of the corresponding rescaled version of the individual-based model, and we compare them with the results of travelling wave analysis (cf. Section~\ref{sec:numerical_results}). We conclude with a discussion of our findings and propose some future research directions (cf. Section~\ref{sec:conclusion}).

\section{The individual-based model}
\label{sec:model}
In this section, integrating the modelling approaches that we developed in~\cite{bubba2020discrete,macfarlane2022individual}, we formulate a phenotype-structured individual-based haptotaxis model of cancer invasion. In this model, individual cells are represented as agents, while the density of ECM and the concentration of MDEs are described by non-negative functions. We allow cells to undergo undirected random movement, phenotype-dependent haptotactic movement in response to the ECM, heritable spontaneous phenotypic changes (i.e. heritable phenotypic changes that occur randomly and are not biased by the cell microenvironment), and phenotype-density-dependent proliferation (i.e. division and death). We consider the scenario where cells also perform phenotype-dependent secretion of MDEs, which then diffuse throughout the spatial domain according to Fick's first Law of diffusion and undergo natural decay. Furthermore, the MDEs break down the ECM to create a gradient that affects the haptotaxis of cancer cells. 

Focussing on a 1D spatial scenario, we let the cells, the density of ECM, and the concentration of MDEs be distributed along the real line $\mathbb{R}$. Furthermore, we describe the phenotypic state of each cell by means of a structuring variable $y\in[0,Y]\subset\mathbb{R}^+$, which takes into account the intercellular variability in haptotactic sensitivity, MDE secretion rate, and proliferation rate. 

\begin{figure}[h!]
\centering
\includegraphics[width=0.4\textwidth]{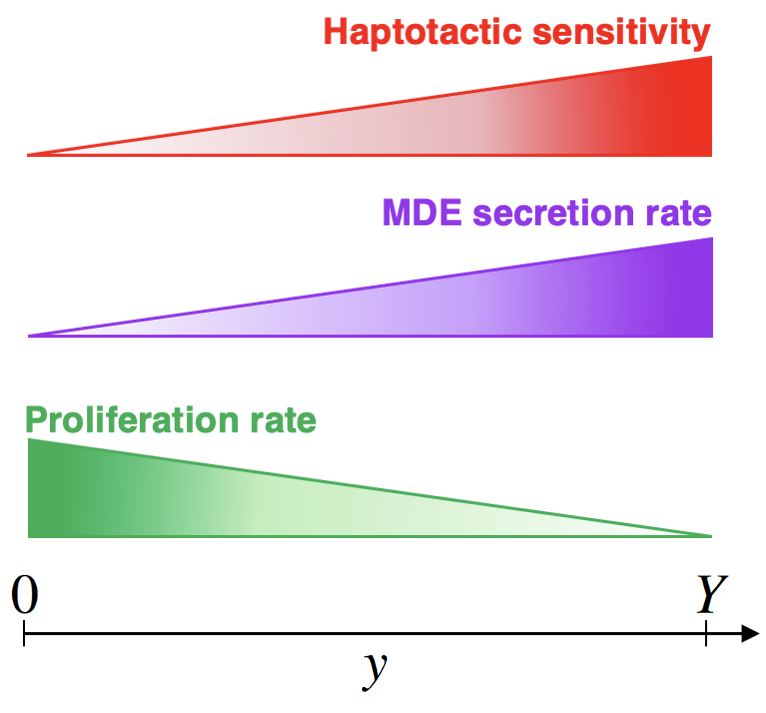}
\caption{{ Schematic overview of the trade-offs between haptotactic and MDE production ability and proliferative potential incorporated into the individual-based model.} }
\label{fig:schpheno}
\end{figure}

In particular, we consider the case where larger values of the structuring variable $y$ correspond to a higher ability to sense spatial gradients of ECM and produce MDEs but a lower proliferative potential (cf. Figure~\ref{fig:schpheno}). This choice is motivated by the energetic costs associated with enhanced motility and greater capacity to alter the surrounding microenvironment, which lead to trade-offs in the form of reduced proliferative potential~\cite{konen2017image}.
 
Hence, cells in phenotypic states characterised by values of $y$ closer to $0$ display a more epithelial-like phenotype (i.e. they behave more like follower cells), whereas cells in phenotypic states characterised by values of $y$ closer to $Y$ display a more mesenchymal-like phenotype (i.e. they behave more like leader cells).

We discretise the time variable $t\in\mathbb{R}^+$ and the space variable $x\in\mathbb{R}$, respectively, as $t_k=k\uptau$ and $x_i=i\Delta_x$ with $k\in \mathbb{N}_0$, $\uptau \in \mathbb{R}^+_*$, $i\in\mathbb{Z}$, and $\Delta_x \in \mathbb{R}^+_*$, where $\mathbb{R}^+_*$ denotes the set of positive real numbers. Moreover, we discretise the phenotype variable via $y_j=j\Delta_y\in[0,Y]$ with $j\in\mathbb{N}_0$ and $\Delta_y \in \mathbb{R}^+_*$. Here, $\uptau$, $\Delta_x$, and $\Delta_y$ are the time-step, space-step, and phenotype-step, respectively.

Each individual cell is represented as an agent that occupies a position on the lattice $\{x_i\}_{i\in\mathbb{Z}}\times\{y_j\}_{j\in\mathbb{N}_0}$ and we introduce the dependent variable $N_{i,j}^k\in\mathbb{N}_0$ to model the number of cells in the phenotypic state $y_j$ at position $x_i$ at time $t_k$. The cell population density and the corresponding cell density are then defined, respectively, as
\begin{equation}\label{eq:nNIB}
n_{i,j}^k\equiv n(t_k,x_i,y_j):=\frac{N_{i,j}^k}{\Delta_x \Delta_y} \quad \text{and} \quad \rho_i^k\equiv\rho(t_k,x_i):=\Delta_y \sum_j n_{i,j}^k.
\end{equation}
The concentration of MDEs and the density of ECM at position $x_i$ at time $t_k$ are denoted by $M_i^k \equiv M(t_k,x_i)$ and $E_i^k \equiv E(t_k,x_i)$, respectively.

\begin{figure}[h!]
\centering
\includegraphics[width=1\textwidth]{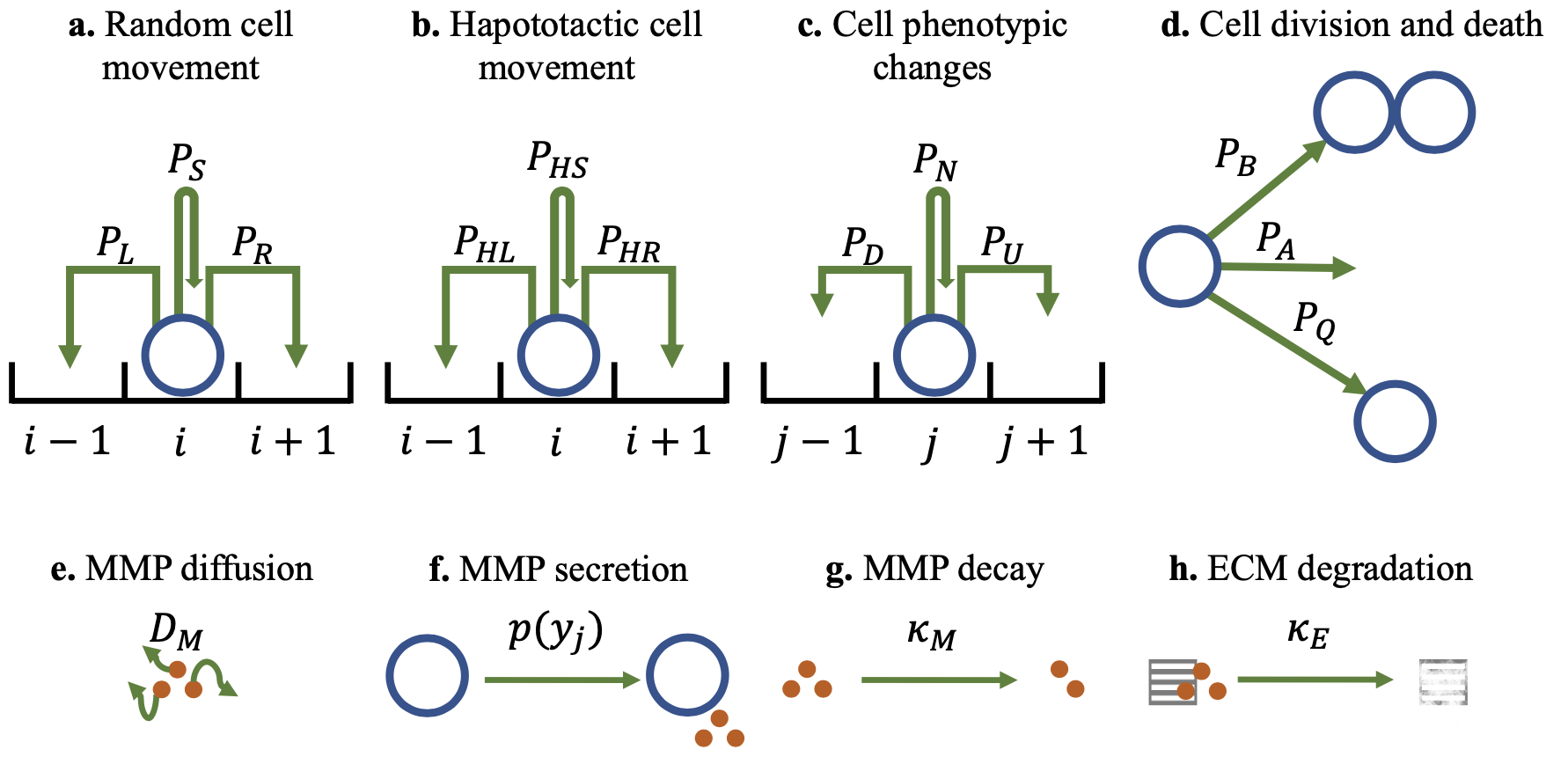}

\caption{{Schematic overview of the mechanisms incorporated into the individual-based model along with the corresponding modelling strategies.} Between time-steps $k$ and $k+1$, each cell in phenotypic state $y_j\in(0,Y)$ at spatial position $x_i\in\mathbb{R}$ may: {\bf a.} move via random motion to either of the positions $x_{i-1}$ and $x_{i+1}$ with probabilities ${P}_{L_{i,j}}^k$ and ${P}_{R_{i,j}}^k$ defined via~\eqref{eq:define_randmov} {or do not undergo random movement with probability ${P}_{S_{i,j}}^k=1-\left({P}_{L_{i,j}}^k+{P}_{R_{i,j}}^k\right)$}; {\bf b.} move via haptotactic motion to either of the positions $x_{i-1}$ and $x_{i+1}$ with probabilities ${P}_{HL_{i,j}}^k$ and ${P}_{HR_{i,j}}^k$ defined via~\eqref{eq:define_haptomov} {or do not undergo haptotactic movement with probability ${P}_{HS_{i,j}}^k=1-\left({P}_{HL_{i,j}}^k+{P}_{HR_{i,j}}^k\right)$}; {\bf c.} undergo a phenotypic change and thus enter either of the phenotype states $y_{j-1}$ and $y_{j+1}$ with probabilities ${P}_{D_{i,j}}^k$ and ${P}_{U_{i,j}}^k$ defined via~\eqref{eq:define_phenochange} {or remain in the same phenotypic state with probability ${P}_{N_{i,j}}^k=1-\left({P}_{D_{i,j}}^k+{P}_{U_{i,j}}^k\right)$}; {\bf d.} die and divide with probabilities ${P}_{A_{i,j}}^k$ and ${P}_{B_{i,j}}^k$ defined via~\eqref{eq:define_proliferation} {or remain quiescent with probability ${P}_{Q_{i,j}}^k=1-\left({P}_{A_{i,j}}^k+{P}_{B_{i,j}}^k\right)$}. The concentration of MDEs will change over time through: {\bf e.} diffusion at the rate $D_M$; {\bf f.} secretion by cells at the phenotypic-dependent rate $p$; {\bf g.} natural decay at rate $\kappa_M$. The ECM density will change over time through: {\bf h.} degradation by MDEs at rate $\kappa_E$.}
\label{fig:scheme}
\end{figure}

The biological mechanisms incorporated into the model and the corresponding modelling strategies are summarised by the schematics in Figure~\ref{fig:scheme} and are described in the remainder of this section.

\subsection{Modelling the dynamics of cancer cells}
\label{sec:cell_model}
As summarised in Figures~\ref{fig:scheme}{\bf a.}-{\bf d.}, between time-steps $k$ and $k+1$, each cell in phenotypic state $y_j\in(0,Y)$ at position $x_i \in \mathbb{R}$ can undergo undirected random movement and haptotactic movement (which are regarded as independent processes), heritable spontaneous phenotypic changes, and cell division and death according to the rules provided in the following subsections.

\subsubsection{Random cell movement}
We model undirected cell movement as a random walk along the spatial dimension, with movement probability $0<\theta\le 1$. In particular, for a focal cell in the phenotypic state $y_j$ at spatial position $x_i$ at time $t_k$, we define the probability of moving left or right to spatial positions $x_{i-1}$ or $x_{i+1}$ as ${P}_{L_{i,j}}^k$ or ${P}_{R_{i,j}}^k$, respectively. As we consider this random movement to be undirected and not affected by the cell phenotype, we define
\begin{equation}
{P}_{L_{i,j}}^k={P}_{R_{i,j}}^k :=\frac{\theta}{2}.
\label{eq:define_randmov}
\end{equation}
Note that cells will not undergo random movement with probability 
\[
P_{S_{i,j}}^k := 1-\left({P}_{L_{i,j}}^k+{P}_{R_{i,j}}^k\right).
\]

\subsubsection{Haptotactic cell movement}
We model haptotactic cell movement in response to the ECM as a biased random walk along the spatial dimension. Since cells move up the gradient of the ECM (i.e. they move towards higher ECM densities), we let the haptotactic movement probabilities depend on the difference between the ECM density at the position occupied by the cell and the ECM density at neighbouring positions. Furthermore, we consider the case where larger values of $y_j$ correlate with a higher haptotaxis sensitivity (cf. Figure~\ref{fig:schpheno}). Hence, we modulate the probabilities of haptotactic cell movement by the function $\mu(y_j)$, which provides a measure of the sensitivity to matrix adhesivity gradients of cells in phenotypic state $y_j$ and thus satisfies the following assumptions
\begin{equation}
\label{eq:assumptions_mu}
\mu(0)=0, \quad \frac{\mathrm{d} \mu(y)}{\mathrm{d} y}>0 \;\; \text{ for} \quad y\in(0,Y].
\end{equation}
We then assume that between time-steps $k$ and $k+1$ a cell in phenotypic state $y_j$ at position $x_i$ may move to the position $x_{i-1}$ (i.e. move left) with probability ${P}_{HL_{i,j}}^k$ or move to the position $x_{i+1}$ (i.e. move right) with probability ${P}_{HR_{i,j}}^k$, where we define
\begin{equation}
{P}_{HL_{i,j}}^k := \eta \, \mu (y_j) \frac{\left(E^k_{i-1}-E^k_{i}\right)_+}{2 \, E_{\rm max}} , \quad {P}_{HR_{i,j}}^k := \eta \, \mu(y_j) \frac{\left(E^k_{i+1}-E^k_{i}\right)_+}{2 \, E_{\rm max}}, \quad \text{with} \quad (\cdot)_+:=\max(0,\cdot).
\label{eq:define_haptomov}
\end{equation}
Here, $E_{\rm max} \in \mathbb{R}^+_*$ is the maximum value of the ECM density before cell invasion starts (see also Section~\ref{sec:ecm_model}). Moreover, the parameter $\eta\in \mathbb{R}^+_*$ is a scaling factor which we consider small enough to ensure $\eta \, \mu (y_j) \leq 1$. Hence, the quantities defined via~\eqref{eq:define_haptomov} satisfy ${0<{P}_{HL_{i,j}}^k+{P}_{HR_{i,j}}^k \le1}$ for all values of $i$, $j$, and $k$. Note that cells will not undergo haptotactic movement with probability 
\[
P_{HS_{i,j}}^k := 1-\left({P}_{HL_{i,j}}^k+{P}_{HR_{i,j}}^k\right).
\]

\subsubsection{Cell phenotypic changes}
We model phenotypic changes by allowing cells to update their phenotypic states according to a random walk along the phenotypic dimension. Between time-steps $k$ and $k+1$ every cell enters a new phenotypic state with probability $0<\beta\le 1$, or remains in its current phenotypic state with probability $1-\beta$. Since we consider spontaneous phenotypic changes, we assume that a cell originally in phenotypic state $y_j$ enters state $y_{j-1}$ with probability ${P}_{D_{i,j}}^k$ or enters state $y_{j+1}$ with probability ${P}_{U_{i,j}}^k$, where we define
\begin{equation}
{P}_{D_{i,j}}^k={P}_{U_{i,j}}^k :=\frac{\beta}{2}.
\label{eq:define_phenochange}
\end{equation}
Therefore cells will not undergo phenotypic changes with probability 
\[
{{P}_{N_{i,j}}^k:= 1-\left({P}_{D_{i,j}}^k+{P}_{U_{i,j}}^k\right)}.
\]
Moreover, no-flux boundary conditions are implemented by aborting any attempted phenotypic change of a cell if it requires moving into a phenotypic state outside of the interval $[0,Y]$.

\subsubsection{Cell division and death}
To incorporate the effects of cell proliferation, we assume that a dividing cell is instantly replaced by two identical progeny cells that inherit the spatial position and phenotypic state of the parent cell. Conversely, a cell undergoing cell death is instantly removed from the population. To take into account phenotypic heterogeneity along with density-dependent inhibition of growth, at time-step $t_k$ we assume that the probabilities of division and death for a cell at spatial position $x_i$ depend both on the phenotypic state of the cell and the local cell density $\rho_i^k$.

In particular, to define the probabilities of cell division and death, we introduce the function $R(y_j,\rho_i^k)$, which describes the net growth rate of the cell population density at spatial position $x_i$ and time $t_k$ due to division and death of cells in the phenotypic state $y_j$, and assume that between time-steps $k$ and $k+1$ a cell in phenotypic state $y_j$ at position $x_i$ may die with probability ${P}_{A_{i,j}}^k$, divide with probability ${P}_{B_{i,j}}^k$, or remain quiescent with probability ${P}_{Q_{i,j}}^k := 1-\left({P}_{A_{i,j}}^k+{P}_{B_{i,j}}^k\right)$, where 
\begin{equation}
\label{eq:define_proliferation}
{P}_{A_{i,j}}^k:=\uptau\ R(y_j,\rho_i^k)_-, \quad {P}_{B_{i,j}}^k :=\uptau\ R(y_j,\rho_i^k)_+, \quad \text{with} \quad (\cdot)_-:=-\min(0,\cdot), \quad (\cdot)_+:=\max(0,\cdot).
\end{equation}
By considering the time-step $\uptau$ sufficiently small, we ensure ${P}_{A_{i,j}}^k+{P}_{B_{i,j}}^k\le 1$ for all values of $i$, $j$, and $k$

We consider the scenario where: larger values of $y_j$ correlate with a lower cell proliferation rate (cf. Figure~\ref{fig:schpheno}); cells stop dividing if the cell density at their current position becomes larger than a critical value $\rho_{\rm max} \in \mathbb{R}^+_*$. Therefore, we make the following assumptions
\begin{equation}
\label{eq:assumptions_R}
R(Y,0)=0,\quad R(0,\rho_{\rm max})=0,\quad \frac{\partial R(y,\rho)}{\partial \rho}<0 \quad \text{and}\quad \frac{\partial R(y,\rho)}{\partial y}<0\quad \text{for } (y,\rho)\in(0,Y)\times\mathbb{R}^+.
\end{equation} 

In particular, we focus on a similar case to that considered in~\cite{macfarlane2022individual}, that is, we assume
\begin{equation}
\label{eq:define_R}
R(y,\rho):=\alpha\left(r(y)-\frac{\rho}{\rho_{\rm max}}\right), \quad r(Y)=0,\quad r(0)=1, \quad \text{and} \quad \frac{\mathrm{d}r(y)}{\mathrm{d}y}<0\quad \text{for } y\in(0,Y),
\end{equation}
with $\alpha \in \mathbb{R}^+_*$. Note that, under assumptions \eqref{eq:assumptions_R}, the definitions given by~\eqref{eq:define_proliferation} ensure that if $\rho_i^k\ge \rho_{\rm max}$ then every cell at position $x_i$ can only die or remain quiescent between time-steps $k$ and $k+1$. Hence, throughout the rest of the paper we will assume 
\begin{equation}
\label{eq:maxden_IBIC}
\max_{i\in\mathbb{Z}}\ \rho_i^0\le \rho_{\rm max} 
\end{equation}
so that
\begin{equation}
\label{eq:maxden_IB}
\rho_i^k\le\rho_{\rm max} \quad \text{for all }(k,i)\in\mathbb{N}_0\times\mathbb{Z}.
\end{equation}

\subsection{Modelling the dynamics of the MDEs and the ECM}
\label{sec:mde_ecm_model}
The dynamics of the MDE concentration and the ECM density are governed by the rules provided in the the following subsections, which are summarised by the schematics in Figures~\ref{fig:scheme}{\bf e.}-{\bf h.} and are coupled with the individual-based model for the dynamics of cancer cells that is presented in Section~\ref{sec:cell_model}. 

\subsubsection{Dynamics of the MDE concentration}
We let $D_M \in \mathbb{R}^+_*$ be the diffusivity of the MDEs and we denote by $\kappa_M \in \mathbb{R}^+_*$ the rate at which the MDEs undergo natural decay. To incorporate into the model the secretion of MDEs by cells in the phenotypic state $y_j$, we introduce the function $p(y_j)$. We focus on the scenario where larger values of $y_j$ correlate with a higher MDE secretion rate (cf. Figure~\ref{fig:schpheno}), i.e. we make the assumptions
\begin{equation}
\label{eq:assumptions_p}
p(0)=p_{\rm min} \in \mathbb{R}^+_*, \quad \frac{ \mathrm{d}p(y)}{\mathrm{d}y}>0 \quad \text{for }y\in(0,Y).
\end{equation}
In this framework, the principle of mass balance gives us the following difference equation for the concentration of MDEs
\begin{equation}
\label{eq:MDE}
M_i^{k+1}=M_{i}^k + \uptau\left[ D_M\ (\mathcal{L}\ M^k)_i- \kappa_M \ M_i^k +\Delta_y \sum_{j} \left(p(y_j)\ n_{i,j}^k\right)\right],
\end{equation}
where $\mathcal{L}$ is the finite-difference Laplacian on the lattice $\{x_i\}_{i\in\mathbb{Z}}$, that is,
\[ (\mathcal{L}\ M^k)_i:=\frac{M_{i+1}^k+M_{i-1}^k-2\ M_i^k}{\Delta_x^2}.\]

\subsubsection{Dynamics of the ECM density}
\label{sec:ecm_model}
We denote by $\kappa_E \in \mathbb{R}^+_*$ the rate at which the ECM is degraded by MDEs. The principle of mass balance gives us the following difference equation for the density of ECM
\begin{equation}
\label{eq:ECM}
E_i^{k+1}=E_i^k-\uptau\ \kappa_E \ M_i^k \, E_i^k.
\end{equation}
Recalling that, as mentioned earlier, $E_{\rm max} \in \mathbb{R}^+_*$ is the maximum value of the ECM density before cell invasion starts, we complement the difference equation~\eqref{eq:ECM} with an initial condition such that
\begin{equation}
\label{eq:ECM_IBIC}
\max_{i\in\mathbb{Z}}\ E_i^0\le E_{\rm max} 
\end{equation}
so that
\begin{equation}
\label{eq:maxECM_IBIC}
E_i^k \le E_{\rm max} \quad \text{for all }(k,i)\in\mathbb{N}_0\times\mathbb{Z}.
\end{equation}

\section{The corresponding continuum model}\label{sec:contmodel}
Through an extension of the limiting procedure that we previously employed in~\cite{bubba2020discrete,chaplain2020bridging,macfarlane2020hybrid,macfarlane2022individual}, letting the time-step $\uptau\rightarrow 0$, the space-step $\Delta_x\rightarrow 0$, and the phenotype-step $\Delta_y\rightarrow0$ in such a way that
\begin{equation}\label{eq:scaledparam}
\frac{\Delta_x^2 \theta}{2 \uptau}\rightarrow D \in \mathbb{R}^+_{*}, \quad \frac{\Delta_x^2 \eta}{2 E_{\rm max} \uptau}\rightarrow \nu \in \mathbb{R}^+_*,\quad \text{and} \quad\frac{\Delta_y^2 \beta}{2\uptau}\rightarrow \lambda \in \mathbb{R}^+_*,
\end{equation}
and introducing the definition 
\begin{equation}
\label{eq:defchimu}
\chi(y) := \nu \, \mu(y),
\end{equation}
one can formally show (cf. Appendix~\ref{app:derivation}) that the deterministic continuum counterpart of the individual-based model presented in Section~\ref{sec:model} is given by the PIDE~\eqref{eq:PDE}$_1$ for the local cell population density function, $n(t,x,y)$, subject to zero Neumann (i.e. no-flux) boundary conditions at $y=0$ and $y=Y$, complemented with the relation~\eqref{eq:PDE}$_2$ for the cell density, $\rho(t,x)$, and coupled with the PDE~\eqref{eq:PDE}$_3$ for the MDE concentration, $M(t,x)$, along with the infinite-dimensional ODE~\eqref{eq:PDE}$_4$ for the ECM density, $E(t,x)$. Consistently with assumptions~\eqref{eq:maxden_IBIC} and~\eqref{eq:ECM_IBIC}, the PIDE~\eqref{eq:PDE}$_1$ and the infinite-dimensional ODE~\eqref{eq:PDE}$_4$ are subject to some initial conditions such that
\begin{equation}
\label{eq:IC_PDEs}
\max_{x \in \mathbb{R}} \int_0^Y n(0,x,y) \ \mathrm{d}y \le \rho_{\rm max}, \quad \max_{x \in \mathbb{R}} E(0,x) \le E_{\rm max}.
\end{equation}

\section{Formal asymptotic analysis}
\label{sec:twanalysis}
In this section, building on the formal asymptotic method that we developed in~\cite{lorenzi2022trade,lorenzi2021}, which relies on the Hamilton-Jacobi approach developed in~\cite{barles2009concentration,diekmann2005dynamics,lorz2011dirac,perthame2008dirac}, we carry out travelling wave analysis of an appropriately rescaled version of the model~\eqref{eq:PDE}.

\subsection{Rescaled model}
We focus on a biological scenario wherein cell proliferation, cell production of MDEs, and ECM degradation have a stronger impact on the dynamics of the system than haptotactic cell movement and diffusion of MDEs, which in turn have a stronger impact than random cell movement and cell phenotypic changes~\cite{anderson1998continuous,anderson2000mathematical,huang2013genetic,orsolits2021new,textor2013analytical}. To this end, we introduce a small parameter $\varepsilon \in \mathbb{R}^+_*$ and choose the parameter scaling
\begin{equation}
\label{eq:scaledpde_param}
\nu := \e, \quad D_M := \e, \quad D := \e^2, \quad \lambda := \e^2.
\end{equation}
Moreover, in order to explore the long-time behaviour of the system, we use the time scaling $\displaystyle{t \to \frac{t}{\varepsilon}}$ in the model~\eqref{eq:PDE}. In so doing, recalling the definition given by~\eqref{eq:defchimu}, we obtain the following rescaled system for the local cell population density function, $n_{\varepsilon}(t,x,y)$, the MDE concentration, $M_{\varepsilon}(t,x)$, and the ECM density, $E_{\varepsilon}(t,x)$:
\begin{equation}
\label{eq:PDE_scaled}
\begin{cases}
\displaystyle{\e \, \partial_t n_{\varepsilon} - \e \, \partial_x\Big(\e \, \partial_x n_{\varepsilon} - n_{\varepsilon} \, \mu(y) \, \partial_x E_{\varepsilon}\Big) = R(y,\rho_{\varepsilon}) \, n_{\varepsilon} + \e^2 \, \partial_{yy} n_{\varepsilon}},\\ \ \\
\displaystyle{\rho_{\varepsilon}(t,x) := \int_0^Y n_{\varepsilon}(t,x,y)\ \mathrm{d}y},\\ \ \\
\displaystyle{\e \, \partial_t M_{\varepsilon} - \e \, \partial_{xx}M_{\varepsilon} = \int_0^Y p(y) \, n_{\varepsilon}(t,x,y) \ \mathrm{d}y- \kappa_M \, M_{\varepsilon}},\\ \ \\
\displaystyle{\e \, \partial_t E_{\varepsilon}=-\kappa_E \, E_{\varepsilon} \, M_{\varepsilon}},
\end{cases}
\quad
(x,y) \in \mathbb{R} \times (0,Y).
\end{equation}

\subsection{Formal limit for $\e \to 0$}
We make the real phase WKB ansatz~\cite{barles1989wavefront,evans1989pde,fleming1986pde}
\begin{equation} \label{WKB}
n_{\varepsilon}(t,x,y) = e^{\frac{u_{\varepsilon}(t,x,y)}{\varepsilon}},
\end{equation}
which gives
\[
\partial_t n_{\varepsilon} = \frac{\partial_t u_{\varepsilon}}{\varepsilon} n_{\varepsilon}, \quad \partial_x n_{\varepsilon} = \frac{\partial_x u_{\varepsilon}}{\varepsilon} n_{\varepsilon}, \quad \partial^2_{yy} n_{\varepsilon} = \left(\frac{1}{\varepsilon^2} \left(\partial_y u_{\varepsilon} \right)^2 + \frac{1}{\varepsilon} \partial^2_{yy} u_{\varepsilon} \right) n_{\varepsilon}.
\]
Substituting the above expressions into the PIDE~\eqref{eq:PDE_scaled}$_1$ for $n_\varepsilon$ and rearranging terms gives the following Hamilton-Jacobi equation for $u_{\varepsilon} \equiv u_{\varepsilon}(t,x,y)$ 
\begin{eqnarray*}
\partial_t u_{\varepsilon} + \mu(y) \, \partial_{x} E_{\e} \, \partial_{x} u_{\e} &=& R(y,\rho_{\varepsilon}) + \left(\partial_x u_{\varepsilon} \right)^2 + \left(\partial_y u_{\varepsilon} \right)^2 + \\
&& + \, \varepsilon \left(\partial^2_{xx} u_{\varepsilon} - \mu(y) \, \partial^2_{xx} E_{\e} + \partial^2_{yy} u_{\varepsilon} \right), \quad (x,y) \in \mathbb{R} \times(0,Y).
\end{eqnarray*}
Now let $\rho(t,x)$ be the leading-order term of the asymptotic expansion for $\rho_{\varepsilon}(t,x)$ as $\varepsilon \to 0$. Considering $x \in \mathbb{R}$ such that $\rho(t,x) >0$ (i.e. $x \in {\rm Supp}(\rho)$) and letting $\varepsilon \to 0$ in the above PDE we formally obtain the following equation for the leading-order term $u \equiv u(t,x,y)$ of the asymptotic expansion for $u_{\varepsilon}(t,x,y)$
\begin{equation}
\label{eq:PDEu}
\partial_t u + \mu(y) \, \partial_{x} E \, \partial_{x} u = R(y,\rho) + \left(\partial_x u \right)^2 + \left(\partial_y u \right)^2, \quad (x,y) \in {{\rm Supp}(\rho)} \times(0,Y),
\end{equation}
where $E \equiv E(t,x)$ is the leading-order term of the asymptotic expansion for $E_{\varepsilon}(t,x)$. 

\paragraph{{Constraint on $u$}}
When $\rho_{\e} < \infty$ for all $\e \in \mathbb{R}^+_*$, if $u_{\e}$ is a strictly concave function of $y$ and $u$ is also a strictly concave function of $y$ whose unique maximum point is $\bar{y}(t,x)$ then considering $x \in {\rm Supp}(\rho)$ and letting $\varepsilon \to 0$ in~\eqref{WKB} formally gives the following constraint on $u$
\begin{equation}
\label{eq:ubaryiszero}
u(t,x,\bar{y}(t,x)) = \max_{y \in [0,Y]} u(t,x,y) = 0, \quad x \in {\rm Supp}(\rho),
\end{equation}
which implies that
\begin{equation}
\label{eq:uybaryiszero}
\partial_y u(t,x,\bar{y}(t,x)) = 0 \quad \text{and} \quad \partial_x u(t,x,\bar{y}(t,x)) = 0, \quad x \in {\rm Supp}(\rho).
\end{equation}

\paragraph{{Relation between $\bar{y}(t,x)$ and $\rho(t,x)$}} Evaluating~\eqref{eq:PDEu} at $y=\bar{y}(t,x)$ and using~\eqref{eq:ubaryiszero} and~\eqref{eq:uybaryiszero} we find 
\[
R(\bar{y}(t,x),\rho(t,x)) = 0, \quad x \in {\rm Supp}(\rho),
\]
from which, using the fact that the function $R(y,\rho)$ is defined via~\eqref{eq:define_R}, we obtain the following formula
\begin{equation}
\label{eq:Riszero}
\rho(t,x) = \rho_{\rm max} \, r(\bar{y}(t,x)), \quad x \in {\rm Supp}(\rho).
\end{equation}
The monotonicity assumption~\eqref{eq:define_R} ensures that~\eqref{eq:Riszero} gives a one-to-one correspondence between $\bar{y}(t,x)$ and $\rho(t,x)$.

\paragraph{{Expressions of $M(t,x)$ and $E(t,x)$}} When $n_{\varepsilon}$ is in the form~\eqref{WKB}, if $u_{\e}$ is a strictly concave function of $y$ and $u$ is also a strictly concave function of $y$ that satisfies the constraint~\eqref{eq:ubaryiszero} then the following asymptotic result formally holds
\[
n_{\varepsilon}(t,x,y) \xrightharpoonup[\varepsilon \rightarrow 0]{} \rho(t,x) \, \delta_{\bar{y}(t,x)}(y) \quad \text{weakly in measures},
\]
where $\delta_{\bar{y}(t,x)}(y)$ is the Dirac delta centred at $y=\bar{y}(t,x)$. In this case, focussing on a biological scenario wherein the ECM density is at the maximum level $E_{\rm max}$ before cell invasion starts at $t=0$, letting $\varepsilon \to 0$ in the PDE~\eqref{eq:PDE_scaled}$_3$ for $M_\varepsilon$ and in the infinite-dimensional ODE~\eqref{eq:PDE_scaled}$_4$ for $E_\varepsilon$ we formally obtain the following expressions of the leading-order terms of the asymptotic expansions for $M_{\varepsilon}(t,x)$ and $E_{\varepsilon}(t,x)$
\begin{equation}
\label{eq:lotME}
M(t,x) = \dfrac{p(\bar{y}(t,x))}{\kappa_M} \rho(t,x), \quad E(t,x) = E_{\rm max} \left(1 - \mathbbm{1}_{{\rm Supp}(M(t,\cdot))}(x)\right),
\end{equation}
where $\mathbbm{1}_{(\cdot)}$ denotes the indicator function of the set $(\cdot)$. 

{
\begin{remark}
Note that the behaviour of $E(t,x)$ depicted by~\eqref{eq:lotME} shares similarities with the behaviour of the nutrient concentration in the model analysed in~\cite{jabin2023collective}.
\end{remark}
}

\paragraph{{Transport equation for $\bar{y}$}}
When $R(y,\rho)$ is defined via~\eqref{eq:define_R}, differentiating~\eqref{eq:PDEu} with respect to $y$, evaluating the resulting equation at $y=\bar{y}(t,x)$, and using~\eqref{eq:ubaryiszero} and~\eqref{eq:uybaryiszero} yields
\begin{equation}
\label{eq:PDEuatbary}
\partial^2_{yt} u(t,x,\bar{y}) + \mu(\bar{y}) \, \partial_{x} E \, \partial^2_{yx} u(t,x,\bar{y}) = \partial_{y} r(\bar{y}), \quad x\in {\rm Supp}(\rho).
\end{equation}
Moreover, differentiating~\eqref{eq:uybaryiszero} with respect to $t$ and $x$ we find, respectively,
\[
\partial^2_{ty} u(t,x,\bar{y}) + \partial^2_{yy} u(t,x,\bar{y}) \, \partial_{t} \bar{y}(t,x) = 0 \; \Rightarrow \; \partial^2_{yt} u(t,x,\bar{y}) = - \partial^2_{yy} u(t,x,\bar{y}) \, \partial_{t} \bar{y}(t,x), \quad x\in {\rm Supp}(\rho)
\]
and
\[
\partial^2_{xy} u(t,x,\bar{y}) + \partial^2_{yy} u(t,x,\bar{y}) \, \partial_{x} \bar{y}(t,x) = 0 \; \Rightarrow \; \partial^2_{yx} u(t,x,\bar{y}) = - \partial^2_{yy} u(t,x,\bar{y}) \, \partial_{x} \bar{y}(t,x), \quad x\in {\rm Supp}(\rho).
\]
Substituting the above expressions of $\partial^2_{yt} u(t,x,\bar{y})$ and $\partial^2_{yx} u(t,x,\bar{y})$ into~\eqref{eq:PDEuatbary}, and using the fact that if $u$ is a strictly concave function of $y$ whose unique maximum point is $\bar{y}(t,x)$ then $\partial^2_{yy} u(t,x,\bar{y}) < 0$, gives the following transport equation for $\bar{y}(t,x)$ 
\begin{equation}
\label{eq:PDEbary}
\partial_{t} \bar{y} + \mu(\bar{y}) \, \partial_{x} E \, \partial_{x} \bar{y} = \frac{\partial_{y} r(\bar{y})}{-\partial^2_{yy} u(t,x,\bar{y})}, \quad x \in {\rm Supp}(\rho).
\end{equation}

\subsection{Travelling wave analysis}
\paragraph{{Travelling wave problem}} 
Substituting the travelling-wave ansatz
\[
u(t,x,y) = u(z,y), \quad \rho(t,x) = \rho(z), \quad \bar{y}(t,x) = \bar{y}(z), \quad M(t,x) = M(z), \quad \text{and} \quad E(t,x) = E(z)
\]
with
\[
z = x - c \, t, \quad c \in \mathbb{R}^+_*,
\]
into~\eqref{eq:Riszero}, \eqref{eq:lotME}, and~\eqref{eq:PDEbary} gives 
\begin{equation}
\label{eq:TWRiszero}
\rho(z) = \rho_{\rm max} \, r(\bar{y}(z)), \quad z \in {\rm Supp}(\rho),
\end{equation}
\begin{equation}
\label{eq:lotMETW}
M(z) = \dfrac{p(\bar{y}(z))}{\kappa_M} \rho(z), \quad E(z) = E_{\rm max} \left(1 - \mathbbm{1}_{{\rm Supp}(M)}(z)\right),
\end{equation}
\begin{equation}
\label{eq:TWbaryori}
\left(c - \mu(\bar{y}) E' \right) \bar{y}' = \frac{\partial_{y} r(\bar{y})}{\partial^2_{yy} u(z,\bar{y})}, \quad z \in {\rm Supp}(\rho).
\end{equation}
The expression~\eqref{eq:lotMETW} of $E(z)$ makes it possible to simplify the differential equation~\eqref{eq:TWbaryori} as follows
\begin{equation}
\label{eq:TWbary}
c \, \bar{y}' = \frac{\partial_{y} r(\bar{y})}{\partial^2_{yy} u(z,\bar{y})}, \quad z \in {\rm Supp}(\rho).
\end{equation}

We complement the differential equation~\eqref{eq:TWbary} with the following asymptotic condition
\begin{equation}
\label{eq:TWBCy}
\lim_{z \to - \infty} \bar{y}(z) =0,
\end{equation}
so that, since $r(0)=1$ (cf. assumptions~\eqref{eq:define_R}), the relation~\eqref{eq:TWRiszero} gives 
\begin{equation}
\label{eq:TWBCrho}
\lim_{z \to - \infty} \rho(z) = \rho_{\rm max}.
\end{equation}

\paragraph{{Shape of travelling waves}} 
Since $\partial_y r(y)<0$ for $y \in (0,Y]$ (cf. assumptions~\eqref{eq:define_R}) and given the fact that if $u$ is a strictly concave function of $y$ whose unique maximum point is $\bar{y}(t,x)$ then $\partial^2_{yy} u(t,x,\bar{y}) < 0$, the differential equation~\eqref{eq:TWbary} along with the relation~\eqref{eq:TWRiszero} ensure that
\begin{equation}
\label{eq:TWbaryincrhodec}
\bar{y}'(z) > 0 \quad \text{and} \quad \rho'(z) < 0, \quad z \in {\rm Supp}(\rho).
\end{equation}
The relation~\eqref{eq:TWRiszero} and the monotonicity results~\eqref{eq:TWbaryincrhodec} along with the fact that $r(Y)=0$ (cf. assumptions~\eqref{eq:define_R}) imply that the position of the edge of the travelling front $\rho(z)$ coincides with the unique point $\ell \in \mathbb{R}$ such that $\bar{y}(\ell)=Y$ and $\bar{y}(z) < Y$ on $(-\infty, \ell)$. Hence, ${\rm Supp}(\rho) = (-\infty, \ell)$ and, since $p(y)>0$ for all $y \in [0,Y]$ (cf. assumptions~\eqref{eq:assumptions_p}), the expressions~\eqref{eq:lotMETW} of $M(z)$ and $E(z)$ yield
\begin{equation}
\label{eq:lotMETWred}
M(z) = \rho_{\rm max} \dfrac{p(\bar{y}(z)) \, r(\bar{y}(z))}{\kappa_M} \mathbbm{1}_{(-\infty, \ell)}(z), \quad E(z) = E_{\rm max} \left(1 - \mathbbm{1}_{(-\infty, \ell)}(z)\right).
\end{equation}

\section{Numerical simulations}
\label{sec:numerical_results}
In this section, we report on numerical solutions of the rescaled continuum model~\eqref{eq:PDE_scaled} and numerical simulations of the corresponding rescaled version of the individual-based model, and we compare them with the results of travelling wave analysis presented in the previous section.

\subsection{Set-up of numerical simulations}
\label{sec:numerical_setup}
We start by describing the set-up used to carry out numerical simulations.

\subsubsection{Model functions and parameters}
To allow the individual-based model to represent the same scenario as the rescaled continuum model~\eqref{eq:PDE_scaled}, we use the same time scaling $t_k\ \rightarrow\ \frac{t_k}{\varepsilon}=k\frac{\tau}{\varepsilon}$ and reformulate the governing rules for the cell dynamics detailed in Section~\ref{sec:model} in terms of 
\[
n_{\varepsilon_{i,j}}^k\equiv n_\varepsilon(t_k,x_i,y_j)=n\left(\frac{t_k}{\varepsilon},x_i,y_j\right):=\frac{N_{\varepsilon_{i,j}}^k}{\Delta_x \Delta_y}, \quad \rho_{\varepsilon_{i}}^k\equiv \rho_\varepsilon(t_k,x_i)=\rho\left(\frac{t_k}{\varepsilon},x_i\right):= \Delta_y \sum_j n_{\varepsilon_{i,j}}^k,
\]
\[
M_{\varepsilon_{i}}^k\equiv M_\varepsilon(t_k,x_i)=M\left(\frac{t_k}{\varepsilon},x_i\right),\quad E_{\varepsilon_{i}}^k\equiv E_\varepsilon(t_k,x_i)=E\left(\frac{t_k}{\varepsilon},x_i\right).
\]
To ensure that conditions~\eqref{eq:scaledparam} and~\eqref{eq:scaledpde_param} are simultaneously satisfied, we additionally set
\[
\theta=\frac{2\uptau}{\Delta_x^2}\varepsilon^2, \quad \eta=\frac{2E_{\rm max}\uptau}{\Delta_x^2}\varepsilon, \quad \beta=\frac{2\uptau}{\Delta_y^2}\varepsilon^2.
\]

In order to carry out numerical simulations, we consider the time interval $[0,T]$ with $T=30$. Furthermore, we restrict the physical domain to the interval $[0,X]$, with $X=100$, and choose $Y=1$. Moreover, we specifically choose $\Delta_x=5 \times 10^{-2}$, $\Delta_y=2 \times 10^{-2}$, and $\tau = \dfrac{\Delta_x^2}{2}$.

To satisfy assumptions \eqref{eq:assumptions_mu}, we use the definition
\begin{equation}
\mu(y):=y^2,
\end{equation}
while in order to satisfy assumptions \eqref{eq:assumptions_R} we define $R(y,\rho)$ via~\eqref{eq:define_R} with $\alpha=0.1$ and, having chosen $Y=1$, we further define 
\begin{equation}
r(y) :=1-y^2.
\end{equation}
To satisfy assumptions \eqref{eq:assumptions_p} on $p(y)$ we also define
\begin{equation}
p(y):=p_{\rm min}+ \zeta y^2,
\end{equation} 
where $\zeta=10^{-5}$ and $p_{\rm min}=10^{-7}$. Furthermore, in the simulations we choose $\kappa_M=1$, $\kappa_E=1$, and $E_{\rm{max}}=1$.

\subsubsection{Initial conditions}
We consider a biological scenario in which, initially, the cell population is localised along the $x=0$ boundary of the spatial domain and most of the cells are in the phenotypic state $y=\bar{y}^0$ at every position $x \in [0,X]$. Specifically, we implement the following initial cell distribution for the IB model 
\begin{equation}
\label{def:icIB}
N_{i,j}^0=\text{int}(F(x_i,y_j)) \; \text{ with } \; F(x,y):= A_0 \ C\ e^{-x^2} \ e^{-\frac{\left(y-\bar{y}^0\right)^2}{\varepsilon}}, 
\end{equation}
where $\text{int}(\cdot)$ is the integer part of $(\cdot)$ and $C$ is a normalisation constant such that 
\[
C\int_0^Y e^{-\frac{\left(y-\bar{y}^0\right)^2}{\varepsilon}} \ \mathrm{d}y=1.
\]
We choose $\bar{y}^0=0.2$ and $A_0=100$. The initial cell density $\rho_{i}^0$ is then calculated from~\eqref{def:icIB} according to the definition given by~\eqref{eq:nNIB}, and we set $\displaystyle{\rho_{\rm max}=\max_{i} \rho^0_i}$.

Moreover, we assume that there are initially no MDEs and the density of ECM is uniform, that is,
\[
M^0_i=0 \quad \text{and} \quad E^0_i=E_{\rm{max}}\quad \text{for all } i.
\]
Finally, we consider different values of $\e$, that is, $\e \in \left\{10^{-2}, 5 \times 10^{-3}, 10^{-3}\right\}$, in order to verify whether, for $\e$ small enough, there is a good agreement between the results of numerical simulations and the results of formal asymptotic analysis for $\e \to 0$ presented in Section~\ref{sec:twanalysis}.

\subsection{Computational implementation of the individual-based model}
All simulations of the individual-based model were performed in {\sc{Matlab}} and the numerical scheme used to solve the rescaled system~\eqref{eq:PDE_scaled} was also implemented in {\sc{Matlab}}.

In the individual-based model, at each time-step, every single cell undergoes a four-step process: (i) random cell movement, according to the probabilities defined via~\eqref{eq:define_randmov}; (ii) haptotactic cell movement, according to the probabilities defined via~\eqref{eq:define_haptomov}; (iii) phenotypic changes, according to the probabilities defined via~\eqref{eq:define_phenochange}; (iv) cell division and death, according to the probabilities defined via~\eqref{eq:define_proliferation}. For every single cell, during each step of this process, a random number is drawn from the standard uniform distribution on the interval $(0,1)$ using the built-in {\sc{Matlab}} function {\sc{rand}}. It is then evaluated whether this number is lower than the probability of the event occurring and if so the event occurs. Since to carry out numerical simulations we have to restrict the spatial domain to the interval $[0, X]$, the attempted movement of a cell is aborted if it requires moving out of this interval. Furthermore, the concentration of MDEs and the density of ECM are calculated using the discrete difference equations~\eqref{eq:MDE} and~\eqref{eq:ECM}, respectively.

\subsection{Numerical methods for the continuum model}
To solve numerically the rescaled system~\eqref{eq:PDE_scaled} posed on $(0, T] \times (0,X) \times (0,Y)$ subject to zero-flux boundary conditions and complemented with the continuum analogues of the initial conditions selected for the individual-based model, we employ a finite volume scheme modified from our previous work~\cite{bubba2020discrete}. 

\subsection{Main results of numerical simulations}
\begin{figure}[h!]
\includegraphics[width=\textwidth]{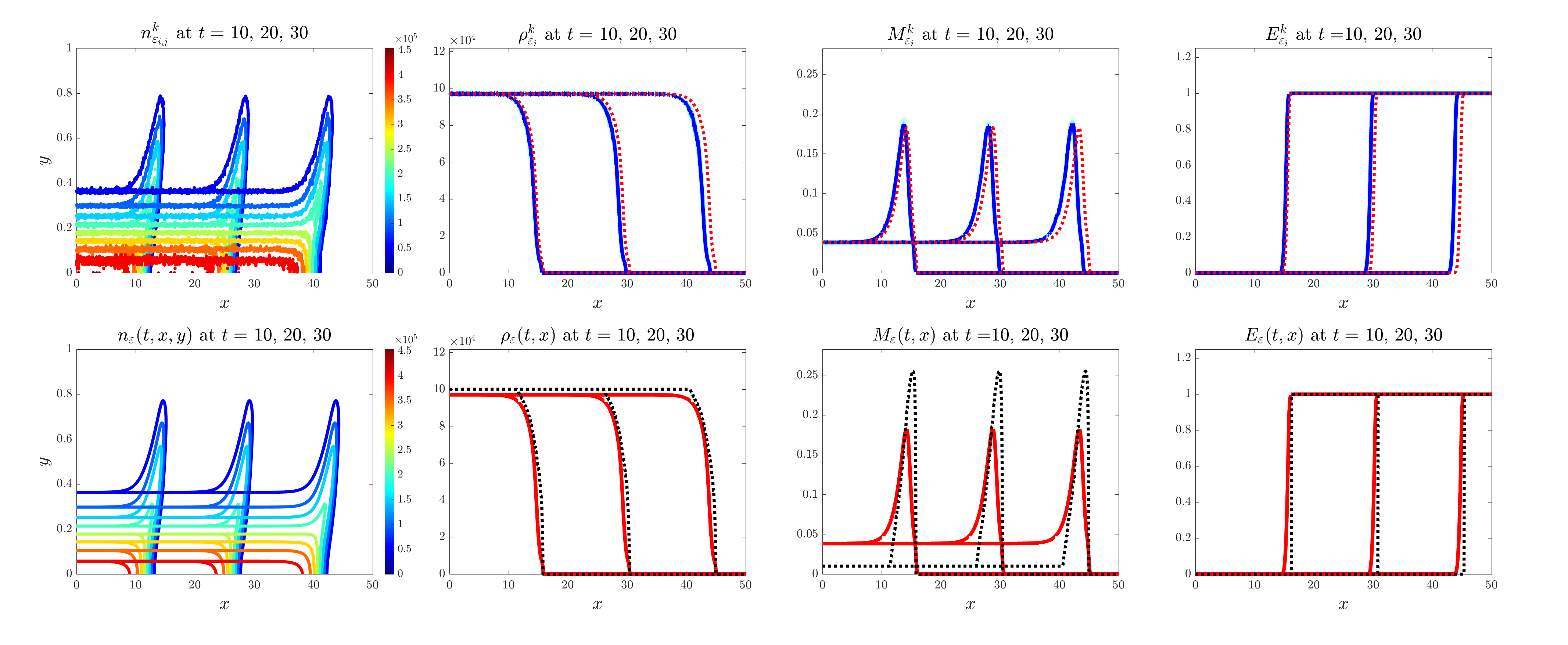}
\caption{Numerical simulation results of the individual-based model (top row) and numerical solutions of the corresponding continuum model~\eqref{eq:PDE_scaled} (bottom row) in the case where $\varepsilon=10^{-2}$. The plots display, from left to right, the cell population density, $n_\varepsilon$, the cell density, $\rho_\varepsilon$, the MMP concentration, $M_\varepsilon$, and the ECM density, $E_\varepsilon$, at progressive times (i.e. $t = 10$, $t = 20$, and $t = 30$) for both modelling approaches. {\bf{ Top row.}} The results from the individual-based model were obtained by averaging over 5 simulations (solid blue lines), and we additionally plot the corresponding results of each simulation (solid cyan lines). We also include the equivalent numerical solutions of the continuum model (dotted red lines) for comparison. {\bf{ Bottom row.}} The values of $\rho_\varepsilon$, $M_\varepsilon$, and $E_\varepsilon$ (solid red lines) are plotted along with the corresponding values obtained by substituting $\bar{y}(t,x)=\bar{y}_{\e}(t,x)$ into the formulas given by~\eqref{eq:Riszero} and~\eqref{eq:lotME} (black dotted lines), with $\bar{y}_{\e}(t,x)$ being the maximum point of the numerical solution $n_\varepsilon(t,x,y)$ to the PIDE~\eqref{eq:PDE_scaled}$_1$ at position $x$ at time $t$.}
\label{fig:resultsA3}
\end{figure}

\begin{figure}[h!]
\includegraphics[width=\textwidth]{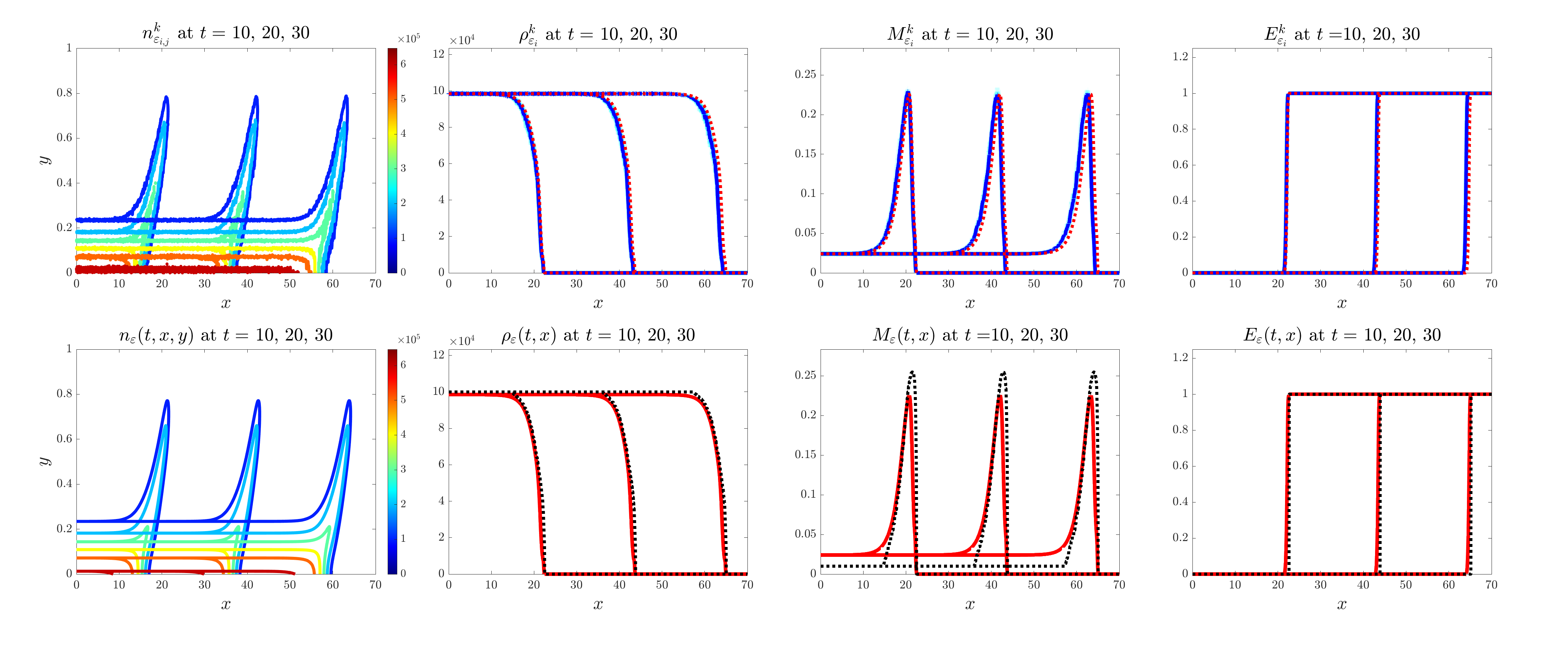}
\caption{Numerical simulation results of the individual-based model (top row) and numerical solutions of the corresponding continuum model~\eqref{eq:PDE_scaled} (bottom row) in the case where $\varepsilon=5 \times 10^{-3}$. The plots display, from left to right, the cell population density, $n_\varepsilon$, the cell density, $\rho_\varepsilon$, the MMP concentration, $M_\varepsilon$, and the ECM density, $E_\varepsilon$, at progressive times (i.e. $t = 10$, $t = 20$, and $t = 30$) for both modelling approaches. {\bf{ Top row.}} The results from the individual-based model were obtained by averaging over 5 simulations (solid blue lines), and we additionally plot the corresponding results of each simulation (solid cyan lines). We also include the equivalent numerical solutions of the continuum model (dotted red lines) for comparison. {\bf{ Bottom row.}} The values of $\rho_\varepsilon$, $M_\varepsilon$, and $E_\varepsilon$ (solid red lines) are plotted along with the corresponding values obtained by substituting $\bar{y}(t,x)=\bar{y}_{\e}(t,x)$ into the formulas given by~\eqref{eq:Riszero} and~\eqref{eq:lotME} (black dotted lines), with $\bar{y}_{\e}(t,x)$ being the maximum point of the numerical solution $n_\varepsilon(t,x,y)$ to the PIDE~\eqref{eq:PDE_scaled}$_1$ at position $x$ at time $t$.}
\label{fig:resultsB3}
\end{figure}

\begin{figure}[h!]
\includegraphics[width=\textwidth]{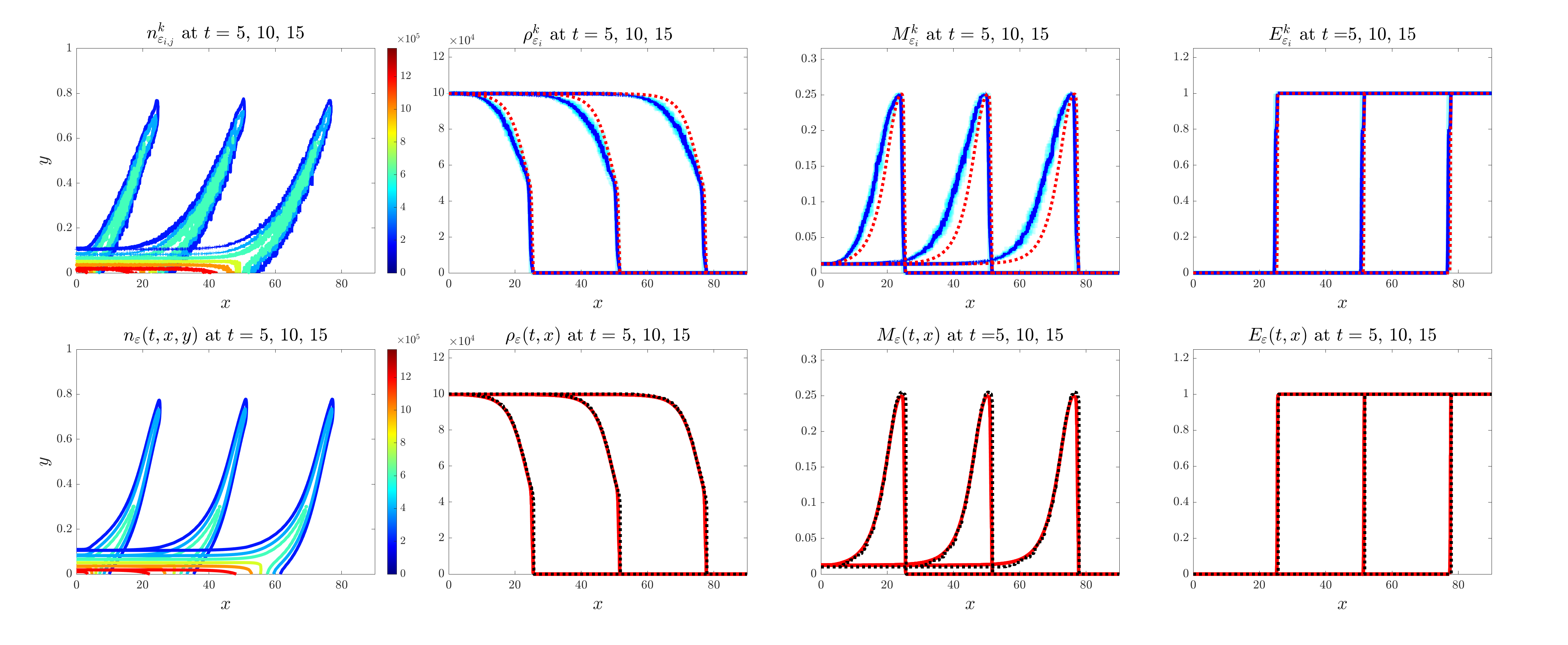}
\caption{Numerical simulation results of the individual-based model (top row) and numerical solutions of the corresponding continuum model~\eqref{eq:PDE_scaled} (bottom row) in the case where $\varepsilon=10^{-3}$. The plots display, from left to right, the cell population density, $n_\varepsilon$, the cell density, $\rho_\varepsilon$, the MMP concentration, $M_\varepsilon$, and the ECM density, $E_\varepsilon$, at progressive times (i.e. $t = 5$, $t = 10$, and $t = 15$) for both modelling approaches. {\bf{ Top row.}} The results from the individual-based model were obtained by averaging over 5 simulations (solid blue lines), and we additionally plot the corresponding results of each simulation (solid cyan lines). We also include the equivalent numerical solutions of the continuum model (dotted red lines) for comparison. {\bf{ Bottom row.}} The values of $\rho_\varepsilon$, $M_\varepsilon$, and $E_\varepsilon$ (solid red lines) are plotted along with the corresponding values obtained by substituting $\bar{y}(t,x)=\bar{y}_{\e}(t,x)$ into the formulas given by~\eqref{eq:Riszero} and~\eqref{eq:lotME} (black dotted lines), with $\bar{y}_{\e}(t,x)$ being the maximum point of the numerical solution $n_\varepsilon(t,x,y)$ to the PIDE~\eqref{eq:PDE_scaled}$_1$ at position $x$ at time $t$.}
\label{fig:resultsC3}
\end{figure}

Our main results of the numerical simulations of the individual-based model and the corresponding numerical solutions of the continuum model~\eqref{eq:PDE_scaled} for three distinct values of the scaling parameter $\varepsilon$ are summarised by the plots in Figures~\ref{fig:resultsA3}-\ref{fig:resultsC3}, which correspond to $\varepsilon=10^{-2}$, $\varepsilon=5 \times 10^{-3}$, and $\varepsilon=10^{-3}$, respectively.

The plots in the top rows of Figures~\ref{fig:resultsA3}-\ref{fig:resultsC3} are the results of the individual-based model averaged over 5 simulations. In particular, from left to right, we have the cell population density, $n_{\varepsilon_{i,j}}^k$, the cell density, $\rho_{\varepsilon_{i}}^k$, the MDE concentration, $M_{\varepsilon_{i}}^k$, and the ECM density, $E_{\varepsilon_{i}}^k$, at progressive times. On the other hand, the plots in the bottom rows of Figures~\ref{fig:resultsA3}-\ref{fig:resultsC3} are the numerical solutions of the continuum model, which are plotted along with the corresponding analytical results presented in Section~\ref{sec:twanalysis}.

These plots show that there is a good agreement between the results of numerical simulations of the individual-based model and numerical solutions of the continuum model. This validates the limiting procedure that we employed to formally derive the continuum model. The same plots also demonstrate that the smaller is the value of $\e$, then the better the agreement between numerical solutions of the rescaled continuum model and the analytical results presented in Section~\ref{sec:twanalysis}. This validates the formal asymptotic method that we used to construct, in the limit as $\e \to 0$, invading fronts with spatial structuring of cell phenotypes. In particular, the plots in Figure~\ref{fig:resultsC3} demonstrate that, when $\e$ is sufficiently small: 
\begin{itemize}
\item[(i)] The local cell population density function $n_{\e}(t,x,y)$ becomes concentrated as a sharp Gaussian with maximum at a point $\bar{y}_{\e}(t,x)$ for all $x$ where $\rho_{\e}(t,x)>0$.
\item[(ii)] The maximum point $\bar{y}_{\e}(t,x)$ behaves like a compactly supported and monotonically increasing travelling front that connects $0$ to $Y$ -- recall that here $Y=1$. This indicates that cells in phenotypic states $y\approx Y$ are concentrated towards the leading edge of the invading front, while cells in phenotypic states corresponding to smaller values of $y$ make up the bulk of the population in the rear.
\item[(iii)] The cell density $\rho_{\e}(t,x)$ behaves like a one-sided compactly supported and monotonically decreasing travelling front that connects $\rho_{\rm max}$ to $0$.
\item[(iv)] The values of $\rho_\varepsilon$, $M_\varepsilon$, and $E_\varepsilon$ are consistent with the values obtained by substituting $\bar{y}(t,x)=\bar{y}_{\e}(t,x)$ into the formulas given by~\eqref{eq:Riszero} and~\eqref{eq:lotME}.
\end{itemize}

\newpage
\section{Conclusion}
\label{sec:conclusion}
We have formulated a model for cancer invasion in which the infiltrating cancer cells can occupy a spectrum of states in the phenotype space, ranging from `fully mesenchymal’ to `fully epithelial’. More precisely, the more mesenchymal cells are those that display stronger haptotaxis responses and have greater capacity to modify the ECM through enhanced secretion of MDEs. However, as a trade-off, they have lower proliferative capacity than the more epithelial cells. The framework is multiscale in that we start with an individual-based model that tracks the dynamics of single cells and is based on a branching random walk over a lattice, where cell movements take place through both physical and phenotype space. By applying limiting techniques, we have formally derived the corresponding continuum model, which takes the form of system~\eqref{eq:PDE}. Despite the intricacy of the model, we showed, through formal asymptotic techniques, that for certain parameter regimes it is possible to carry out a detailed travelling wave analysis and obtain the form for the wave profile. Simulations have been performed of both the individual-based model and the continuum model, which generally show excellent correspondence. Moreover, when parameter values are chosen from the appropriate parameter regime, numerical solutions to the continuum model match closely with the corresponding analytical form. This validates the formal limiting procedure that we employed to derive the continuum model alongside the formal asymptotic method that we used to characterise the wave profile.

Notably, solutions to the model reveal a capacity for self-organisation, in the sense that an initially almost homogeneous population resolves itself into an invading front with spatial structuring of phenotypes. Precisely, the most mesenchymal cells dominate the leading edge of the invasion wave and the most epithelial (and most proliferative) dominate the rear, representing a bulk tumour population. As such, the model recapitulates similar observations into a front to back structuring of invasion waves into leader-type and follower-type cells, witnessed in an increasing number of experimental studies over recent years \cite{vilchez2021decoding}.

A number of other continuum models have been formulated to study how phenotypic diversity alters cancer invasion processes. These include those intended to describe ``go-or-grow’’ dynamics, a term coined for glioma growth processes where a dichotomy of cells into proliferating or migrating cell types has been suggested \cite{giese1996dichotomy}. In many of these models, heterogeneity is restricted to binary states (a proliferating class and a migrating class), with functions defining the state to state transitions. Often these models have restricted to relatively simple assumptions for cell migration processes, e.g. a simple diffusion process \cite{pham2012density,stepien2018traveling}, although more complex movement models have also been considered, e.g. \cite{conte2021}. The model here expands the potential framework for modelling go-or-growth processes, to cover all potential states between fully proliferative and fully migratory. Binary cell state models with distinct proliferative and migratory characteristics have also been developed to describe acid-mediated cancer invasion \cite{strobl2020mix}, where a similar wave structuring of the distinct populations can be found under certain configurations. Invasion models with continuous phenotypes that range from more migratory to more proliferative states have been applied to avascular cancer growth, where further variables are included to describe tissue oxygen levels \cite{fiandaca2022phenotype}. The study here provides a structure for more detailed analytical studies into the travelling wave dynamics observed in these primarily numerically-based investigations. While not specifically focussing on cancer invasion, continuous phenotypic structuring models have also been developed and analysed for chemotaxis-driven wave invasion in~\cite{lorenzi2022trade} and density- or pressure-driven wave invasion in~\cite{lorenzi2021,macfarlane2022individual}; the study here expands on the methodologies introduced therein, reinforcing their utility to study a diverse range of models used to explain invasion processes. 

\begin{figure}
\includegraphics[width=\textwidth]{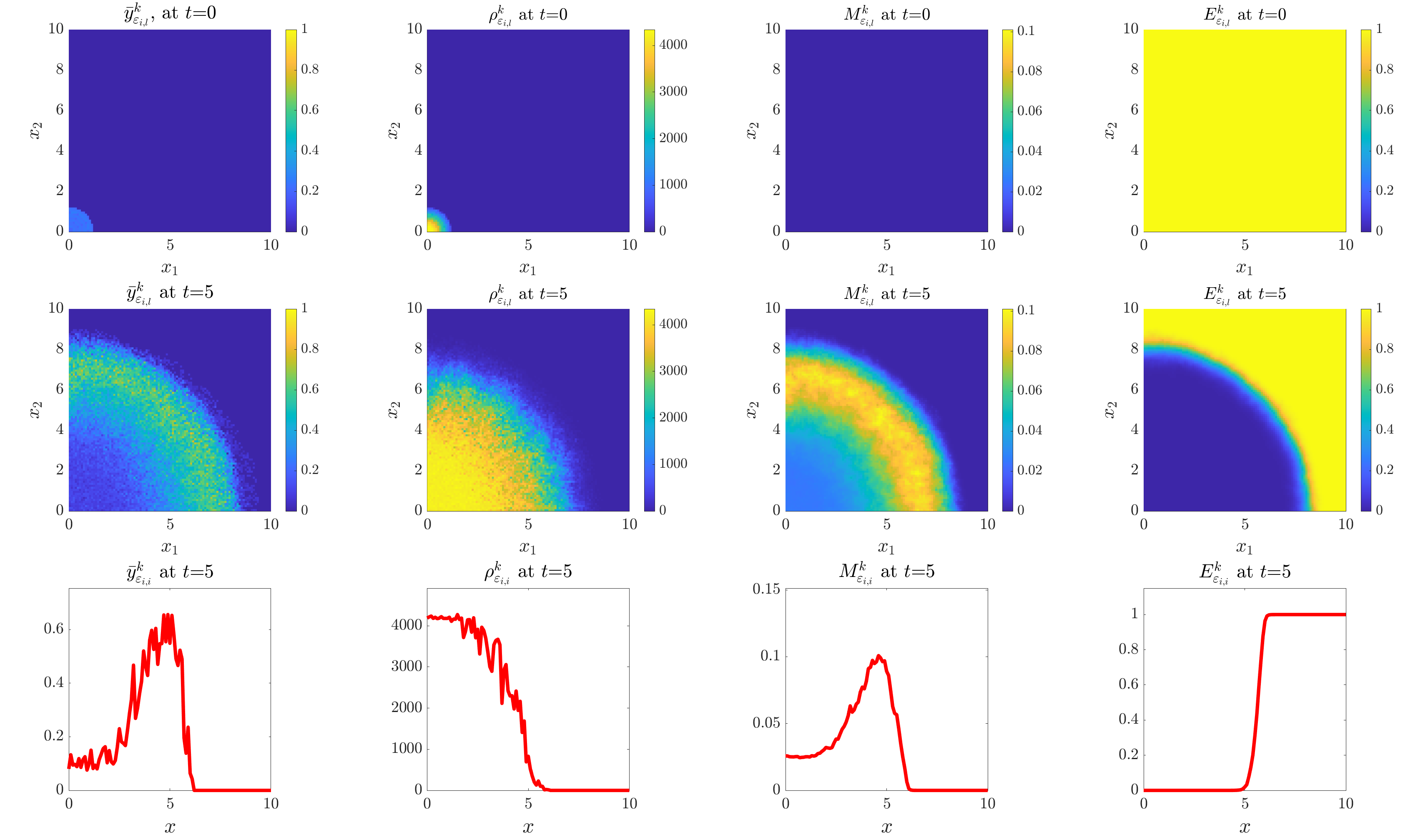}
\caption{Preliminary 2D results from the individual-based model, obtained by averaging over 15 simulations, in the case where $\varepsilon=10^{-2}$. {\bf{Top and middle rows.}} The plots display, from left to right, the maximum point of the cell population density, ${\bar{y}_{\e} = \argmax_{y \in [0,Y]} n_\varepsilon}$, the cell density, $\rho_\varepsilon$, the MMP concentration, $M_\varepsilon$, and the ECM density, $E_\varepsilon$, at the start and end of simulations -- i.e. $t = 0$ (top) and $t = 5$ (middle). {\bf{Bottom row.}} The plots display, from left to right, $\displaystyle{\bar{y}_{\e}}$, $\rho_\varepsilon$, $M_\varepsilon$, and $E_\varepsilon$ across the radial transect at the end of simulations (i.e. at  $t = 5$). Here $x\equiv(x_1,x_2) \in [0,10] \times [0,10]$ with $\Delta_{x_1}=\Delta_{x_2}=0.1$, $y \in [0,1]$ with $\Delta_y=0.02$, and the initial cell distribution is the 2D analogue of~\eqref{def:icIB} with $A_0=1$, while all the other parameters and functions are kept the same as in the 1D simulations of Figures~\ref{fig:resultsA3}-\ref{fig:resultsC3}.}
\label{fig:results2D}
\end{figure}

There are, clearly, various further extensions that could be considered. The current study has concentrated on invasion processes in 1D for each of physical space (i.e. a transect across the invading front) and phenotype space (from epithelial-like to mesenchymal-like). It is possible, of course, to extend the dimensionality of either or both of these spaces. As a way of illustration, preliminary simulations are presented in Figure \ref{fig:results2D} for an extension to 2D for the physical space in the individual-based model, where from an initial population concentrated at the origin we observe the emergence of a quasi-symmetric growing tumour with leader-to-follower structuring across the radial transect (consistent with the corresponding 1D model). A natural question to explore in two dimensions would be whether there are conditions under which the symmetric growth breaks, e.g. whether spatial structuring such as `tumour fingering’ emerges. Previous models have shown that fingering can occur in various scenarios, such as tumour infiltration into heterogeneous ECM environments~\cite{anderson2006tumor} and tumour growth in the presence of cells with different mobilities~\cite{drasdo2012modeling,lorenzi2017interfaces}. Whether it is also possible for such phenomena to develop under certain assumptions for the manner in which phenotypic transitions occur could be a point of focus. 

Extensions in the dimensionality of the phenotype space may also be of interest. 
Here we have considered a linear pathway from the fully epithelial state to the fully mesenchymal state -- i.e. where reduced proliferation is accompanied simultaneously by an upregulation in both haptotaxis and secretion of MDEs. Given that EMTs in cancers can often be `partial’ in nature \cite{vilchez2021decoding}, it is possible that different pathways could be taken from epithelial to mesenchymal: cells could follow separate pathways in which first haptotactic movement is upregulated and then secretion of MDEs, or vice versa. Extensions to a higher dimensionality in the phenotype space would allow exploration into whether this can give rise to additional subtlety in the positioning of different phenotypes.

Another natural extension would be to target the model towards particular experimental studies of leader-follower behaviour, by incorporating system-specific phenomena. For example, studies of collective invasion in non-small cell lung cancer (NSCLC) tumour spheroids have led to an experimental model of invasion with intricate signalling between leader and follower cells \cite{konen2017image}: in terms of movement, leader cells secrete fibronectin and release VEGF that guides follower cell movements through chemotaxis; in terms of proliferation, followers secrete factors to promote leader cell proliferation, while leaders secrete factors that hinder follower growth. Adapting the model to include these additional factors and their impact on follower/leader behaviour would provide a means to test the experimental model and, for example, investigate how perturbations to various aspects of the signalling system would impact on the rate of infiltration. 

Further exploration could be made into the processes that lead to phenotypic changes. Here, we have adopted the relatively simple assumption of (unbiased) random phenotype switches, which at the cell-population level leads to diffusion across phenotype space. It is also quite plausible that phenotype changes may be biased in particular directions -- e.g. from epithelial towards mesenchymal (or vice versa) -- and that the direction and strength of the bias changes with the tumour microenvironment \cite{aggarwal2021interplay}. Our model has also decoupled proliferation from phenotypic changes, effectively assuming that proliferative events lead to daughter cells of the same phenotype; divisions may also occur in an asymmetric manner -- e.g. division of a follower cell leading to a daughter of leader phenotype. With our framework, the impact of such changes can be investigated both at the individual and continuous level.

Summarising, the work here provides the framework for developing and analysing sophisticated haptotaxis models for cancer invasion in which the cell population contains significant phenotypic heterogeneity. While these models present significant challenges at both a numerical and analytical level, we believe the methods developed and described here can allow for further progress in this area. 

\section*{Acknowledgments}
TL gratefully acknowledges support from the Italian Ministry of University and Research (MUR) through the grant PRIN 2020 project (No. 2020JLWP23) ``Integrated Mathematical Approaches to Socio-Epidemiological Dynamics'' (CUP: E15F21005420006) and the grant PRIN2022-PNRR project (No. P2022Z7ZAJ) ``A Unitary Mathematical Framework for Modelling Muscular Dystrophies'' (CUP: E53D23018070001), from the CNRS International Research Project ``Modelisation de la biomecanique cellulaire et tissulaire'' (MOCETIBI), and from the Istituto Nazionale di Alta Matematica (INdAM) and the Gruppo Nazionale per la Fisica Matematica (GNFM). KJP is a member of INdAM-GNFM and acknowledges ``Miur-Dipartimento di Eccellenza'' funding to the Dipartimento di Scienze, Progetto e Politiche del Territorio (DIST). TL and KJP would also like to thank the Isaac Newton Institute for Mathematical Sciences, Cambridge, for support and hospitality during the programme Mathematics of Movement, where work on this paper was undertaken.


\bibliographystyle{ejamike}
\bibliography{references}

\appendix 
\section{Appendices}
{\small
\subsection{Formal derivation of the continuum model~\eqref{eq:PDE}}
\label{app:derivation}
We provide the full details of the formal derivation of the continuum model~\eqref{eq:PDE} from the individual-based model presented in Section~\ref{sec:model}, which relies on an extension of the limiting procedure that we previously employed in~\cite{bubba2020discrete,chaplain2020bridging,macfarlane2020hybrid,macfarlane2022individual}.

\subsubsection{Equation for cell population density}
When cell dynamics are governed by the rules underlying the individual-based model presented in Section~\ref{sec:model}, the principle of mass balance gives, for the cell population density,
\begin{eqnarray}
\label{Der1}
n^{k+1}_{i,j}&=&n^{k}_{i+1,j+1}\left\{ \frac{\beta}{2}\left[ 1+\uptau R(y_j, \rho^{k}_{i})\right]\left[\frac{\theta}{2}+\frac{\eta \mu(y_{j})}{2 E_{\rm max}} \left(E^{k}_{i}-E^{k}_{i+1}\right)_+\right]\right\}\\
&&+n^{k}_{i-1,j+1}\left\{ \frac{\beta}{2}\left[ 1+\uptau R(y_j, \rho^{k}_{i})\right]\left[\frac{\theta}{2}+\frac{\eta \mu(y_{j})}{2 E_{\rm max}} \left(E^{k}_{i}-E^{k}_{i-1}\right)_+\right]\right\}\nonumber\\
&&+n^{k}_{i+1,j-1}\left\{ \frac{\beta}{2}\left[ 1+\uptau R(y_j, \rho^{k}_{i})\right]\left[\frac{\theta}{2}+\frac{\eta \mu(y_{j})}{2 E_{\rm max}} \left(E^{k}_{i}-E^{k}_{i+1}\right)_+\right]\right\}\nonumber\\
&&+n^{k}_{i-1,j-1}\left\{ \frac{\beta}{2}\left[ 1+\uptau R(y_j, \rho^{k}_{i})\right]\left[\frac{\theta}{2}+\frac{\eta \mu(y_{j})}{2 E_{\rm max}} \left(E^{k}_{i}-E^{k}_{i-1}\right)_+\right]\right\}\nonumber\\
&&+n^{k}_{i,j+1}\left\{ \frac{\beta}{2}\left[ 1+\uptau R(y_j, \rho^{k}_{i})\right]\left[1-\theta-\frac{\eta \mu(y_{j})}{2 E_{\rm max}}\left[ \left(E^{k}_{i+1}-E^{k}_{i}\right)_+ + \left(E^{k}_{i-1}-E^{k}_{i}\right)_+\right]\right] \right\}\nonumber\\
&&+n^{k}_{i,j-1}\left\{ \frac{\beta}{2}\left[ 1+\uptau R(y_j, \rho^{k}_{i})\right]\left[1-\theta-\frac{\eta \mu(y_{j})}{2 E_{\rm max}}\left[ \left(E^{k}_{i+1}-E^{k}_{i}\right)_+ + \left(E^{k}_{i-1}-E^{k}_{i}\right)_+\right]\right] \right\}\nonumber\\
&&+n^{k}_{i+1,j}\left\{ (1-\beta)\left[ 1+\uptau R(y_j, \rho^{k}_{i})\right]\left[\frac{\theta}{2}+\frac{\eta \mu(y_{j})}{2 E_{\rm max}} \left(E^{k}_{i}-E^{k}_{i+1}\right)_+\right]\right\}\nonumber\\
&&+n^{k}_{i-1,j}\left\{ (1-\beta)\left[ 1+\uptau R(y_j, \rho^{k}_{i})\right]\left[\frac{\theta}{2}+\frac{\eta \mu(y_{j})}{2 E_{\rm max}} \left(E^{k}_{i}-E^{k}_{i-1}\right)_+\right]\right\}\nonumber\\
&&+n^{k}_{i,j}\left\{ (1-\beta)\left[ 1+\uptau R(y_j, \rho^{k}_{i})\right] \left[1-\theta-\frac{\eta \mu(y_{j})}{2 E_{\rm max}}\left[\left(E^{k}_{i+1}-E^{k}_{i}\right)_++\left(E^{k}_{i-1}-E^{k}_{i}\right)_+ \right]\right]\right\}.\nonumber
\end{eqnarray}
Using the fact that for $\uptau$, $\Delta_x$, and $\Delta_y$ sufficiently small the following relations hold
\begin{eqnarray}
n^{k}_{i,j}\approx n(t,x,y)\equiv n,\quad n^{k+1}_{i,j}\approx n(t+\uptau,x,y),\quad n^{k}_{i\pm1,j}\approx n(t,x\pm\Delta_x,y),\quad n^{k}_{i,j\pm1}\approx n(t,x,y\pm\Delta_y)\nonumber\\
E^{k}_{i}\approx E(t,x)\equiv E, \quad E^{k}_{i\pm1}\approx E(t,x\pm\Delta_x), \quad \mu(y_j)\approx \mu(y) \equiv \mu, \quad \rho^{k}_{i} \approx \rho(t,x)\equiv \rho, \quad R(y_j,\rho^{k}_{i})\approx R(y,\rho) \equiv R,\nonumber
\end{eqnarray}
Eq \eqref{Der1} can be rewritten as
\begin{eqnarray}
\label{Der2}
n(t+\uptau,x,y)&=&n(t,x+\Delta_x,y+\Delta_y)\left\{ \frac{\beta}{2}\left[ 1+\uptau R\right]\left[\frac{\theta}{2}+\frac{\eta \mu}{2 E_{\rm max}} \left(E-E(t,x+\Delta_x)\right)_+\right]\right\}\\
&&+n(t,x-\Delta_x,y+\Delta_y)\left\{ \frac{\beta}{2}\left[ 1+\uptau R\right]\left[\frac{\theta}{2}+\frac{\eta \mu}{2 E_{\rm max}} \left(E-E(t,x-\Delta_x)\right)_+\right]\right\}\nonumber\\
&&+n(t,x+\Delta_x,y-\Delta_y)\left\{ \frac{\beta}{2}\left[ 1+\uptau R\right]\left[\frac{\theta}{2}+\frac{\eta \mu}{2 E_{\rm max}} \left(E-E(t,x+\Delta_x)\right)_+\right]\right\}\nonumber\\
&&+n(t,x-\Delta_x,y-\Delta_y)\left\{ \frac{\beta}{2}\left[ 1+\uptau R\right]\left[\frac{\theta}{2}+\frac{\eta \mu}{2 E_{\rm max}} \left(E-E(t,x-\Delta_x)\right)_+\right]\right\}\nonumber\\
&&+n(t,x,y+\Delta_y)\left\{ \frac{\beta}{2}\left[ 1+\uptau R\right]\left[1-\theta-\frac{\eta \mu}{2 E_{\rm max}}\left[ \left(E(t,x+\Delta_x)-E\right)_+ + \left(E(t,x-\Delta_x)-E\right)_+\right]\right] \right\}\nonumber\\
&&+n(t,x,y-\Delta_y)\left\{ \frac{\beta}{2}\left[ 1+\uptau R\right]\left[1-\theta-\frac{\eta \mu}{2 E_{\rm max}}\left[ \left(E(t,x+\Delta_x)-E\right)_+ + \left(E(t,x-\Delta_x)-E\right)_+\right]\right] \right\}\nonumber\\
&&+n(t,x+\Delta_x,y)\left\{ (1-\beta)\left[ 1+\uptau R\right]\left[\frac{\theta}{2}+\frac{\eta \mu}{2 E_{\rm max}} \left(E-E(t,x+\Delta_x)\right)_+\right]\right\}\nonumber\\
&&+n(t,x-\Delta_x,y)\left\{ (1-\beta)\left[ 1+\uptau R\right]\left[\frac{\theta}{2}+\frac{\eta \mu}{2 E_{\rm max}} \left(E-E(t,x-\Delta_x)\right)_+\right]\right\}\nonumber\\
&&+n\left\{ (1-\beta)\left[ 1+\uptau R\right] \left[1-\theta-\frac{\eta \mu}{2 E_{\rm max}}\left[\left(E(t,x+\Delta_x)-E\right)_++\left(E(t,x-\Delta_x)-E\right)_+ \right]\right]\right\}.\nonumber
\end{eqnarray}

Assuming the function $n$ to be sufficiently regular, we now use the following Taylor expansions
\begin{eqnarray}
n(t,x,y\pm\Delta_y)= n\pm\Delta_y\frac{\partial n}{\partial y}+\frac{\Delta_y^2}{2}\frac{\partial^2 n}{\partial y^2}+h.o.t., \quad n(t,x\pm\Delta_x,y)= n\pm\Delta_x\frac{\partial n}{\partial x}+\frac{\Delta_x^2}{2}\frac{\partial^2 n}{\partial x^2}+h.o.t.,\nonumber\\
n(t,x+\Delta_x,y\pm\Delta_y)= n+\Delta_x\frac{\partial n}{\partial x} \pm\Delta_y\frac{\partial n}{\partial y}+\frac{\Delta_x^2}{2}\frac{\partial^2 n}{\partial x^2}+\frac{\Delta_y^2}{2}\frac{\partial^2 n}{\partial y^2} \pm \Delta_x\Delta_y \frac{\partial^2 n}{\partial x\partial y}+h.o.t. ,\nonumber\\
n(t,x-\Delta_x,y\pm\Delta_y)= n-\Delta_x\frac{\partial n}{\partial x} \pm\Delta_y\frac{\partial n}{\partial y}+\frac{\Delta_x^2}{2}\frac{\partial^2 n}{\partial x^2}+\frac{\Delta_y^2}{2}\frac{\partial^2 n}{\partial y^2} \mp \Delta_x \Delta_y \frac{\partial^2 n}{\partial x\partial y}+h.o.t. ,\nonumber
\end{eqnarray}
which allow us to rewrite \eqref{Der2} as
\begin{eqnarray}
\label{Der3}
&&n(t+\uptau,x,y)=\\
&&n\left\{\frac{\beta}{2}\left[ 1+\uptau R\right]\left[\frac{\theta}{2}+\frac{\eta \mu}{2 E_{\rm max}}\left(E-E(t,x+\Delta_x)\right)_+\right]\right\}+\Delta_x\frac{\partial n}{\partial x}\left\{\frac{\beta}{2}\left[ 1+\uptau R\right]\left[\frac{\theta}{2}+\frac{\eta \mu}{2 E_{\rm max}}\left(E-E(t,x+\Delta_x)\right)_+\right]\right\}\nonumber\\%
&&+\Delta_y\frac{\partial n}{\partial y}\left\{ \frac{\beta}{2}\left[ 1+\uptau R\right]\left[\frac{\theta}{2}+\frac{\eta \mu}{2 E_{\rm max}} \left(E-E(t,x+\Delta_x)\right)_+\right]\right\}\nonumber+\Delta_x \Delta_y \frac{\partial^2 n}{\partial x\partial y}\left\{ \frac{\beta}{2}\left[ 1+\uptau R\right]\left[\frac{\theta}{2}+\frac{\eta \mu}{2 E_{\rm max}} \left(E-E(t,x+\Delta_x)\right)_+\right]\right\}\nonumber\\%
&&+\frac{\Delta_x^2}{2}\frac{\partial^2 n}{\partial x^2}\left\{ \frac{\beta}{2}\left[ 1+\uptau R\right]\left[\frac{\theta}{2}+\frac{\eta \mu}{2 E_{\rm max}} \left(E-E(t,x+\Delta_x)\right)_+\right]\right\}\nonumber+\frac{\Delta_y^2}{2}\frac{\partial^2 n}{\partial y^2}\left\{ \frac{\beta}{2}\left[ 1+\uptau R\right]\left[\frac{\theta}{2}+\frac{\eta \mu}{2 E_{\rm max}} \left(E-E(t,x+\Delta_x)\right)_+\right]\right\}\nonumber\\
&&+n\left\{ \frac{\beta}{2}\left[ 1+\uptau R\right]\left[\frac{\theta}{2}+\frac{\eta \mu}{2 E_{\rm max}} \left(E-E(t,x-\Delta_x)\right)_+\right]\right\}\nonumber-\Delta_x\frac{\partial n}{\partial x}\left\{ \frac{\beta}{2}\left[ 1+\uptau R\right]\left[\frac{\theta}{2}+\frac{\eta \mu}{2 E_{\rm max}} \left(E-E(t,x-\Delta_x)\right)_+\right]\right\}\nonumber\\%
&&+\Delta_y\frac{\partial n}{\partial y}\left\{ \frac{\beta}{2}\left[ 1+\uptau R\right]\left[\frac{\theta}{2}+\frac{\eta \mu}{2 E_{\rm max}} \left(E-E(t,x-\Delta_x)\right)_+\right]\right\}\nonumber-\Delta_x \Delta_y \frac{\partial^2 n}{\partial x\partial y}\left\{ \frac{\beta}{2}\left[ 1+\uptau R\right]\left[\frac{\theta}{2}+\frac{\eta \mu}{2 E_{\rm max}} \left(E-E(t,x-\Delta_x)\right)_+\right]\right\}\nonumber\\%
&&+\frac{\Delta_x^2}{2}\frac{\partial^2 n}{\partial x^2}\left\{ \frac{\beta}{2}\left[ 1+\uptau R\right]\left[\frac{\theta}{2}+\frac{\eta \mu}{2 E_{\rm max}} \left(E-E(t,x-\Delta_x)\right)_+\right]\right\}\nonumber+\frac{\Delta_y^2}{2}\frac{\partial^2 n}{\partial y^2}\left\{ \frac{\beta}{2}\left[ 1+\uptau R\right]\left[\frac{\theta}{2}+\frac{\eta \mu}{2 E_{\rm max}} \left(E-E(t,x-\Delta_x)\right)_+\right]\right\}\nonumber\\%
&&+n\left\{ \frac{\beta}{2}\left[ 1+\uptau R\right]\left[\frac{\theta}{2}+\frac{\eta \mu}{2 E_{\rm max}} \left(E-E(t,x+\Delta_x)\right)_+\right]\right\}\nonumber+\Delta_x\frac{\partial n}{\partial x}\left\{ \frac{\beta}{2}\left[ 1+\uptau R\right]\left[\frac{\theta}{2}+\frac{\eta \mu}{2 E_{\rm max}} \left(E-E(t,x+\Delta_x)\right)_+\right]\right\}\nonumber\\%
&&-\Delta_y\frac{\partial n}{\partial y}\left\{ \frac{\beta}{2}\left[ 1+\uptau R\right]\left[\frac{\theta}{2}+\frac{\eta \mu}{2 E_{\rm max}} \left(E-E(t,x+\Delta_x)\right)_+\right]\right\}\nonumber-\Delta_x \Delta_y \frac{\partial^2 n}{\partial x\partial y}\left\{ \frac{\beta}{2}\left[ 1+\uptau R\right]\left[\frac{\theta}{2}+\frac{\eta \mu}{2 E_{\rm max}} \left(E-E(t,x+\Delta_x)\right)_+\right]\right\}\nonumber\\%
&&+\frac{\Delta_x^2}{2}\frac{\partial^2 n}{\partial x^2}\left\{ \frac{\beta}{2}\left[ 1+\uptau R\right]\left[\frac{\theta}{2}+\frac{\eta \mu}{2 E_{\rm max}} \left(E-E(t,x+\Delta_x)\right)_+\right]\right\}\nonumber+\frac{\Delta_y^2}{2}\frac{\partial^2 n}{\partial y^2}\left\{ \frac{\beta}{2}\left[ 1+\uptau R\right]\left[\frac{\theta}{2}+\frac{\eta \mu}{2 E_{\rm max}} \left(E-E(t,x+\Delta_x)\right)_+\right]\right\}\nonumber\\%
&&+n\left\{ \frac{\beta}{2}\left[ 1+\uptau R\right]\left[\frac{\theta}{2}+\frac{\eta \mu}{2 E_{\rm max}} \left(E-E(t,x-\Delta_x)\right)_+\right]\right\}\nonumber-\Delta_x\frac{\partial n}{\partial x}\left\{ \frac{\beta}{2}\left[ 1+\uptau R\right]\left[\frac{\theta}{2}+\frac{\eta \mu}{2 E_{\rm max}} \left(E-E(t,x-\Delta_x)\right)_+\right]\right\}\nonumber\\%
&&-\Delta_y\frac{\partial n}{\partial y}\left\{ \frac{\beta}{2}\left[ 1+\uptau R\right]\left[\frac{\theta}{2}+\frac{\eta \mu}{2 E_{\rm max}} \left(E-E(t,x-\Delta_x)\right)_+\right]\right\}+\Delta_x \Delta_y \frac{\partial^2 n}{\partial x\partial y}\left\{ \frac{\beta}{2}\left[ 1+\uptau R\right]\left[\frac{\theta}{2}+\frac{\eta \mu}{2 E_{\rm max}} \left(E-E(t,x-\Delta_x)\right)_+\right]\right\}\nonumber\\%
&&+\frac{\Delta_x^2}{2}\frac{\partial^2 n}{\partial x^2}\left\{ \frac{\beta}{2}\left[ 1+\uptau R\right]\left[\frac{\theta}{2}+\frac{\eta \mu}{2 E_{\rm max}} \left(E-E(t,x-\Delta_x)\right)_+\right]\right\}\nonumber+\frac{\Delta_y^2}{2}\frac{\partial^2 n}{\partial y^2}\left\{ \frac{\beta}{2}\left[ 1+\uptau R\right]\left[\frac{\theta}{2}+\frac{\eta \mu}{2 E_{\rm max}} \left(E-E(t,x-\Delta_x)\right)_+\right]\right\}\nonumber\\%
&&+n\left\{ \frac{\beta}{2}\left[ 1+\uptau R\right]\left[1-\theta-\frac{\eta \mu}{2 E_{\rm max}}\left[ \left(E(t,x+\Delta_x)-E\right)_+ + \left(E(t,x-\Delta_x)-E\right)_+\right]\right] \right\}\nonumber\\
&&+\Delta_y\frac{\partial n}{\partial y}\left\{ \frac{\beta}{2}\left[ 1+\uptau R\right]\left[1-\theta-\frac{\eta \mu}{2 E_{\rm max}}\left[ \left(E(t,x+\Delta_x)-E\right)_+ + \left(E(t,x-\Delta_x)-E\right)_+\right]\right] \right\}\nonumber\\
&&+\frac{\Delta_y^2}{2}\frac{\partial^2 n}{\partial y^2}\left\{ \frac{\beta}{2}\left[ 1+\uptau R\right]\left[1-\theta-\frac{\eta \mu}{2 E_{\rm max}}\left[ \left(E(t,x+\Delta_x)-E\right)_+ + \left(E(t,x-\Delta_x)-E\right)_+\right]\right] \right\}\nonumber\\
&&+n\left\{ \frac{\beta}{2}\left[ 1+\uptau R\right]\left[1-\theta-\frac{\eta \mu}{2 E_{\rm max}}\left[ \left(E(t,x+\Delta_x)-E\right)_+ + \left(E(t,x-\Delta_x)-E\right)_+\right]\right] \right\}\nonumber\\
&&-\Delta_y\frac{\partial n}{\partial y}\left\{ \frac{\beta}{2}\left[ 1+\uptau R\right]\left[1-\theta-\frac{\eta \mu}{2 E_{\rm max}}\left[ \left(E(t,x+\Delta_x)-E\right)_+ + \left(E(t,x-\Delta_x)-E\right)_+\right]\right] \right\}\nonumber\\
&&+\frac{\Delta_y^2}{2}\frac{\partial^2 n}{\partial y^2}\left\{ \frac{\beta}{2}\left[ 1+\uptau R\right]\left[1-\theta-\frac{\eta \mu}{2 E_{\rm max}}\left[ \left(E(t,x+\Delta_x)-E\right)_+ + \left(E(t,x-\Delta_x)-E\right)_+\right]\right] \right\}\nonumber\\
&&+n\left\{ (1-\beta)\left[ 1+\uptau R\right]\left[\frac{\theta}{2}+\frac{\eta \mu}{2 E_{\rm max}} \left(E-E(t,x+\Delta_x)\right)_+\right]\right\}\nonumber+\Delta_x\frac{\partial n}{\partial x}\left\{ (1-\beta)\left[ 1+\uptau R\right]\left[\frac{\theta}{2}+\frac{\eta \mu}{2 E_{\rm max}} \left(E-E(t,x+\Delta_x)\right)_+\right]\right\}\nonumber\\%
&&+\frac{\Delta_x^2}{2}\frac{\partial^2 n}{\partial x^2}\left\{ (1-\beta)\left[ 1+\uptau R\right]\left[\frac{\theta}{2}+\frac{\eta \mu}{2 E_{\rm max}} \left(E-E(t,x+\Delta_x)\right)_+\right]\right\}\nonumber\\
&&+n\left\{ (1-\beta)\left[ 1+\uptau R\right]\left[\frac{\theta}{2}+\frac{\eta \mu}{2 E_{\rm max}} \left(E-E(t,x-\Delta_x)\right)_+\right]\right\}-\Delta_x\frac{\partial n}{\partial x}\left\{ (1-\beta)\left[ 1+\uptau R\right]\left[\frac{\theta}{2}+\frac{\eta \mu}{2 E_{\rm max}} \left(E-E(t,x-\Delta_x)\right)_+\right]\right\}\nonumber\\
&&+\frac{\Delta_x^2}{2}\frac{\partial^2 n}{\partial x^2}\left\{ (1-\beta)\left[ 1+\uptau R\right]\left[\frac{\theta}{2}+\frac{\eta \mu}{2 E_{\rm max}} \left(E-E(t,x-\Delta_x)\right)_+\right]\right\}\nonumber\\%
&&+n\left\{ (1-\beta)\left[ 1+\uptau R\right] \left[1-\theta-\frac{\eta \mu}{2 E_{\rm max}}\left[\left(E(t,x+\Delta_x)-E\right)_++\left(E(t,x-\Delta_x)-E\right)_+ \right]\right]\right\} + h.o.t.\nonumber
\end{eqnarray}
Then rearranging and collecting terms of derivatives of $n$, we obtain
\begin{eqnarray}
\label{Der4}
&&n(t+\uptau,x,y)=\\
&&n\left\{ \frac{\beta}{2}\left[ 1+\uptau R\right]\left[\frac{\theta}{2}+\frac{\eta \mu}{2 E_{\rm max}} \left(E-E(t,x+\Delta_x)\right)_+\right]\right\}+n\left\{ \frac{\beta}{2}\left[ 1+\uptau R\right]\left[\frac{\theta}{2}+\frac{\eta \mu}{2 E_{\rm max}} \left(E-E(t,x-\Delta_x)\right)_+\right]\right\}\nonumber\\%
&&+n\left\{ \frac{\beta}{2}\left[ 1+\uptau R\right]\left[\frac{\theta}{2}+\frac{\eta \mu}{2 E_{\rm max}} \left(E-E(t,x+\Delta_x)\right)_+\right]\right\}\nonumber+n\left\{ \frac{\beta}{2}\left[ 1+\uptau R\right]\left[\frac{\theta}{2}+\frac{\eta \mu}{2 E_{\rm max}} \left(E-E(t,x-\Delta_x)\right)_+\right]\right\}\nonumber\\%
&&+n\left\{ \frac{\beta}{2}\left[ 1+\uptau R\right]\left[1-\theta-\frac{\eta \mu}{2 E_{\rm max}}\left[ \left(E(t,x+\Delta_x)-E\right)_+ + \left(E(t,x-\Delta_x)-E\right)_+\right]\right] \right\}\nonumber\\%
&&+n\left\{ \frac{\beta}{2}\left[ 1+\uptau R\right]\left[1-\theta-\frac{\eta \mu}{2 E_{\rm max}}\left[ \left(E(t,x+\Delta_x)-E\right)_+ + \left(E(t,x-\Delta_x)-E\right)_+\right]\right] \right\}\nonumber\\%
&&+n\left\{ (1-\beta)\left[ 1+\uptau R\right]\left[\frac{\theta}{2}+\frac{\eta \mu}{2 E_{\rm max}} \left(E-E(t,x+\Delta_x)\right)_+\right]\right\}\nonumber+n\left\{ (1-\beta)\left[ 1+\uptau R\right]\left[\frac{\theta}{2}+\frac{\eta \mu}{2 E_{\rm max}} \left(E-E(t,x-\Delta_x)\right)_+\right]\right\}\nonumber\\%
&&+n\left\{ (1-\beta)\left[ 1+\uptau R\right] \left[1-\theta-\frac{\eta \mu}{2 E_{\rm max}}\left[\left(E(t,x+\Delta_x)-E\right)_++\left(E(t,x-\Delta_x)-E\right)_+ \right]\right]\right\}\nonumber\\
&&+\Delta_x\frac{\partial n}{\partial x}\left\{ \frac{\beta}{2}\left[ 1+\uptau R\right]\left[\frac{\theta}{2}+\frac{\eta \mu}{2 E_{\rm max}} \left(E-E(t,x+\Delta_x)\right)_+\right]\right\} -\Delta_x\frac{\partial n}{\partial x}\left\{ \frac{\beta}{2}\left[ 1+\uptau R\right]\left[\frac{\theta}{2}+\frac{\eta \mu}{2 E_{\rm max}} \left(E-E(t,x-\Delta_x)\right)_+\right]\right\}\nonumber\\%
&&+\Delta_x\frac{\partial n}{\partial x}\left\{ \frac{\beta}{2}\left[ 1+\uptau R\right]\left[\frac{\theta}{2}+\frac{\eta \mu}{2 E_{\rm max}} \left(E-E(t,x+\Delta_x)\right)_+\right]\right\}\nonumber-\Delta_x\frac{\partial n}{\partial x}\left\{ \frac{\beta}{2}\left[ 1+\uptau R\right]\left[\frac{\theta}{2}+\frac{\eta \mu}{2 E_{\rm max}} \left(E-E(t,x-\Delta_x)\right)_+\right]\right\}\nonumber\\%
&&+\Delta_x\frac{\partial n}{\partial x}\left\{ (1-\beta)\left[ 1+\uptau R\right]\left[\frac{\theta}{2}+\frac{\eta \mu}{2 E_{\rm max}} \left(E-E(t,x+\Delta_x)\right)_+\right]\right\}\nonumber\\
&&-\Delta_x\frac{\partial n}{\partial x}\left\{ (1-\beta)\left[ 1+\uptau R\right]\left[\frac{\theta}{2}+\frac{\eta \mu}{2 E_{\rm max}} \left(E-E(t,x-\Delta_x)\right)_+\right]\right\}\nonumber\\%
&&+\Delta_y\frac{\partial n}{\partial y}\left\{ \frac{\beta}{2}\left[ 1+\uptau R\right]\left[\frac{\theta}{2}+\frac{\eta \mu}{2 E_{\rm max}} \left(E-E(t,x+\Delta_x)\right)_+\right]\right\}\nonumber+\Delta_y\frac{\partial n}{\partial y}\left\{ \frac{\beta}{2}\left[ 1+\uptau R\right]\left[\frac{\theta}{2}+\frac{\eta \mu}{2 E_{\rm max}} \left(E-E(t,x-\Delta_x)\right)_+\right]\right\}\nonumber\\%
&&-\Delta_y\frac{\partial n}{\partial y}\left\{ \frac{\beta}{2}\left[ 1+\uptau R\right]\left[\frac{\theta}{2}+\frac{\eta \mu}{2 E_{\rm max}} \left(E-E(t,x+\Delta_x)\right)_+\right]\right\}-\Delta_y\frac{\partial n}{\partial y}\left\{ \frac{\beta}{2}\left[ 1+\uptau R\right]\left[\frac{\theta}{2}+\frac{\eta \mu}{2 E_{\rm max}} \left(E-E(t,x-\Delta_x)\right)_+\right]\right\}\nonumber\\%
&&+\Delta_y\frac{\partial n}{\partial y}\left\{ \frac{\beta}{2}\left[ 1+\uptau R\right]\left[1-\theta-\frac{\eta \mu}{2 E_{\rm max}}\left[ \left(E(t,x+\Delta_x)-E\right)_+ + \left(E(t,x-\Delta_x)-E\right)_+\right]\right] \right\}\nonumber\\%
&&-\Delta_y\frac{\partial n}{\partial y}\left\{ \frac{\beta}{2}\left[ 1+\uptau R\right]\left[1-\theta-\frac{\eta \mu}{2 E_{\rm max}}\left[ \left(E(t,x+\Delta_x)-E\right)_+ + \left(E(t,x-\Delta_x)-E\right)_+\right]\right] \right\}\nonumber\\%
&&+\Delta_x \Delta_y \frac{\partial^2 n}{\partial x\partial y}\left\{ \frac{\beta}{2}\left[ 1+\uptau R\right]\left[\frac{\theta}{2}+\frac{\eta \mu}{2 E_{\rm max}} \left(E-E(t,x+\Delta_x)\right)_+\right]\right\}\nonumber\\
&&-\Delta_x \Delta_y \frac{\partial^2 n}{\partial x\partial y}\left\{ \frac{\beta}{2}\left[ 1+\uptau R\right]\left[\frac{\theta}{2}+\frac{\eta \mu}{2 E_{\rm max}} \left(E-E(t,x-\Delta_x)\right)_+\right]\right\}\nonumber\\%
&&-\Delta_x \Delta_y \frac{\partial^2 n}{\partial x\partial y}\left\{ \frac{\beta}{2}\left[ 1+\uptau R\right]\left[\frac{\theta}{2}+\frac{\eta \mu}{2 E_{\rm max}} \left(E-E(t,x+\Delta_x)\right)_+\right]\right\}\nonumber\\
&&+\Delta_x \Delta_y \frac{\partial^2 n}{\partial x\partial y}\left\{ \frac{\beta}{2}\left[ 1+\uptau R\right]\left[\frac{\theta}{2}+\frac{\eta \mu}{2 E_{\rm max}} \left(E-E(t,x-\Delta_x)\right)_+\right]\right\}\nonumber\\%
&&+\frac{\Delta_x^2}{2}\frac{\partial^2 n}{\partial x^2}\left\{ \frac{\beta}{2}\left[ 1+\uptau R\right]\left[\frac{\theta}{2}+\frac{\eta \mu}{2 E_{\rm max}} \left(E-E(t,x+\Delta_x)\right)_+\right]\right\}+\frac{\Delta_x^2}{2}\frac{\partial^2 n}{\partial x^2}\left\{ \frac{\beta}{2}\left[ 1+\uptau R\right]\left[\frac{\theta}{2}+\frac{\eta \mu}{2 E_{\rm max}} \left(E-E(t,x-\Delta_x)\right)_+\right]\right\}\nonumber\\%
&&+\frac{\Delta_x^2}{2}\frac{\partial^2 n}{\partial x^2}\left\{ \frac{\beta}{2}\left[ 1+\uptau R\right]\left[\frac{\theta}{2}+\frac{\eta \mu}{2 E_{\rm max}} \left(E-E(t,x+\Delta_x)\right)_+\right]\right\}+\frac{\Delta_x^2}{2}\frac{\partial^2 n}{\partial x^2}\left\{ \frac{\beta}{2}\left[ 1+\uptau R\right]\left[\frac{\theta}{2}+\frac{\eta \mu}{2 E_{\rm max}} \left(E-E(t,x-\Delta_x)\right)_+\right]\right\}\nonumber\\%
&&+\frac{\Delta_x^2}{2}\frac{\partial^2 n}{\partial x^2}\left\{ (1-\beta)\left[ 1+\uptau R\right]\left[\frac{\theta}{2}+\frac{\eta \mu}{2 E_{\rm max}} \left(E-E(t,x+\Delta_x)\right)_+\right]\right\}\nonumber\\
&&+\frac{\Delta_x^2}{2}\frac{\partial^2 n}{\partial x^2}\left\{ (1-\beta)\left[ 1+\uptau R\right]\left[\frac{\theta}{2}+\frac{\eta \mu}{2 E_{\rm max}} \left(E-E(t,x-\Delta_x)\right)_+\right]\right\}\nonumber\\%
&&+\frac{\Delta_y^2}{2}\frac{\partial^2 n}{\partial y^2}\left\{ \frac{\beta}{2}\left[ 1+\uptau R\right]\left[\frac{\theta}{2}+\frac{\eta \mu}{2 E_{\rm max}} \left(E-E(t,x+\Delta_x)\right)_+\right]\right\}+\frac{\Delta_y^2}{2}\frac{\partial^2 n}{\partial y^2}\left\{ \frac{\beta}{2}\left[ 1+\uptau R\right]\left[\frac{\theta}{2}+\frac{\eta \mu}{2 E_{\rm max}} \left(E-E(t,x-\Delta_x)\right)_+\right]\right\}\nonumber\\%
&&+\frac{\Delta_y^2}{2}\frac{\partial^2 n}{\partial y^2}\left\{ \frac{\beta}{2}\left[ 1+\uptau R\right]\left[\frac{\theta}{2}+\frac{\eta \mu}{2 E_{\rm max}} \left(E-E(t,x+\Delta_x)\right)_+\right]\right\}\nonumber+\frac{\Delta_y^2}{2}\frac{\partial^2 n}{\partial y^2}\left\{ \frac{\beta}{2}\left[ 1+\uptau R\right]\left[\frac{\theta}{2}+\frac{\eta \mu}{2 E_{\rm max}} \left(E-E(t,x-\Delta_x)\right)_+\right]\right\}\nonumber\\
&&+\frac{\Delta_y^2}{2}\frac{\partial^2 n}{\partial y^2}\left\{ \frac{\beta}{2}\left[ 1+\uptau R\right]\left[1-\theta-\frac{\eta \mu}{2 E_{\rm max}}\left[ \left(E(t,x+\Delta_x)-E\right)_+ + \left(E(t,x-\Delta_x)-E\right)_+\right]\right] \right\}\nonumber\\
&&+\frac{\Delta_y^2}{2}\frac{\partial^2 n}{\partial y^2}\left\{ \frac{\beta}{2}\left[ 1+\uptau R\right]\left[1-\theta-\frac{\eta \mu}{2 E_{\rm max}}\left[ \left(E(t,x+\Delta_x)-E\right)_+ + \left(E(t,x-\Delta_x)-E\right)_+\right]\right] \right\} +h.o.t.\nonumber
\end{eqnarray}
Further simplifying yields
\begin{eqnarray}
\label{Der5}
&&n(t+\uptau,x,y)=\\&&n\left\{ \beta\left[ 1+\uptau R\right]\left[\frac{\theta}{2}+\frac{\eta \mu}{2 E_{\rm max}} \left(E-E(t,x+\Delta_x)\right)_+\right]\right\}+n\left\{ \beta\left[ 1+\uptau R\right]\left[\frac{\theta}{2}+\frac{\eta \mu}{2 E_{\rm max}} \left(E-E(t,x-\Delta_x)\right)_+\right]\right\}\nonumber\\%
&&+n\left\{ \beta\left[ 1+\uptau R\right]\left[1-\theta-\frac{\eta \mu}{2 E_{\rm max}}\left[ \left(E(t,x+\Delta_x)-E\right)_+ + \left(E(t,x-\Delta_x)-E\right)_+\right]\right] \right\}\nonumber\\
&&+n\left\{ (1-\beta)\left[ 1+\uptau R\right]\left[\frac{\theta}{2}+\frac{\eta \mu}{2 E_{\rm max}} \left(E-E(t,x+\Delta_x)\right)_+\right]\right\}\nonumber+n\left\{ (1-\beta)\left[ 1+\uptau R\right]\left[\frac{\theta}{2}+\frac{\eta \mu}{2 E_{\rm max}} \left(E-E(t,x-\Delta_x)\right)_+\right]\right\}\nonumber\\%
&&+n\left\{ (1-\beta)\left[ 1+\uptau R\right] \left[1-\theta-\frac{\eta \mu}{2 E_{\rm max}}\left[\left(E(t,x+\Delta_x)-E\right)_++\left(E(t,x-\Delta_x)-E\right)_+ \right]\right]\right\}\nonumber\\%
&&+\Delta_x\frac{\partial n}{\partial x}\left\{ \beta\left[ 1+\uptau R\right]\left[\frac{\theta}{2}+\frac{\eta \mu}{2 E_{\rm max}} \left(E-E(t,x+\Delta_x)\right)_+\right]\right\}-\Delta_x\frac{\partial n}{\partial x}\left\{ \beta\left[ 1+\uptau R\right]\left[\frac{\theta}{2}+\frac{\eta \mu}{2 E_{\rm max}} \left(E-E(t,x-\Delta_x)\right)_+\right]\right\}\nonumber\\%
&&+\Delta_x\frac{\partial n}{\partial x}\left\{ (1-\beta)\left[ 1+\uptau R\right]\left[\frac{\theta}{2}+\frac{\eta \mu}{2 E_{\rm max}} \left(E-E(t,x+\Delta_x)\right)_+\right]\right\}\nonumber\\
&&-\Delta_x\frac{\partial n}{\partial x}\left\{ (1-\beta)\left[ 1+\uptau R\right]\left[\frac{\theta}{2}+\frac{\eta \mu}{2 E_{\rm max}} \left(E-E(t,x-\Delta_x)\right)_+\right]\right\}\nonumber\\%
&&+\Delta_y\frac{\partial n}{\partial y}\left\{ \beta\left[ 1+\uptau R\right]\left[\frac{\theta}{2}+\frac{\eta \mu}{2 E_{\rm max}} \left(E-E(t,x+\Delta_x)\right)_+\right]\right\}-\Delta_y\frac{\partial n}{\partial y}\left\{ \beta\left[ 1+\uptau R\right]\left[\frac{\theta}{2}+\frac{\eta \mu}{2 E_{\rm max}} \left(E-E(t,x-\Delta_x)\right)_+\right]\right\}\nonumber\\%
&&+\frac{\Delta_x^2}{2}\frac{\partial^2 n}{\partial x^2}\left\{ \beta\left[ 1+\uptau R\right]\left[\frac{\theta}{2}+\frac{\eta \mu}{2 E_{\rm max}} \left(E-E(t,x+\Delta_x)\right)_+\right]\right\}+\frac{\Delta_x^2}{2}\frac{\partial^2 n}{\partial x^2}\left\{ \beta\left[ 1+\uptau R\right]\left[\frac{\theta}{2}+\frac{\eta \mu}{2 E_{\rm max}} \left(E-E(t,x-\Delta_x)\right)_+\right]\right\}\nonumber\\%
&&+\frac{\Delta_x^2}{2}\frac{\partial^2 n}{\partial x^2}\left\{ (1-\beta)\left[ 1+\uptau R\right]\left[\frac{\theta}{2}+\frac{\eta \mu}{2 E_{\rm max}} \left(E-E(t,x+\Delta_x)\right)_+\right]\right\}\nonumber\\&&+\frac{\Delta_x^2}{2}\frac{\partial^2 n}{\partial x^2}\left\{ (1-\beta)\left[ 1+\uptau R\right]\left[\frac{\theta}{2}+\frac{\eta \mu}{2 E_{\rm max}} \left(E-E(t,x-\Delta_x)\right)_+\right]\right\}\nonumber\\%
&&+\frac{\Delta_y^2}{2}\frac{\partial^2 n}{\partial y^2}\left\{ \beta\left[ 1+\uptau R\right]\left[\frac{\theta}{2}+\frac{\eta \mu}{2 E_{\rm max}} \left(E-E(t,x+\Delta_x)\right)_+\right]\right\}+\frac{\Delta_y^2}{2}\frac{\partial^2 n}{\partial y^2}\left\{ \beta\left[ 1+\uptau R\right]\left[\frac{\theta}{2}+\frac{\eta \mu}{2 E_{\rm max}} \left(E-E(t,x-\Delta_x)\right)_+\right]\right\}\nonumber\\%
&&+\frac{\Delta_y^2}{2}\frac{\partial^2 n}{\partial y^2}\left\{ \beta\left[ 1+\uptau R\right]\left[1-\theta-\frac{\eta \mu}{2 E_{\rm max}}\left[ \left(E(t,x+\Delta_x)-E\right)_+ + \left(E(t,x-\Delta_x)-E\right)_+\right]\right] \right\} +h.o.t.\nonumber
\end{eqnarray}
Cancelling out $\beta$ terms, where possible, we find
\begin{eqnarray}
\label{Der6}
&&n(t+\uptau,x,y)=\\
&&n\left\{ \left[ 1+\uptau R\right]\left[\frac{\theta}{2}+\frac{\eta \mu}{2 E_{\rm max}} \left(E-E(t,x+\Delta_x)\right)_+\right]\right\}+n\left\{ \left[ 1+\uptau R\right]\left[\frac{\theta}{2}+\frac{\eta \mu}{2 E_{\rm max}} \left(E-E(t,x-\Delta_x)\right)_+\right]\right\}\nonumber\\%
&&+n\left\{ \left[ 1+\uptau R\right]\left[1-\theta-\frac{\eta \mu}{2 E_{\rm max}}\left[ \left(E(t,x+\Delta_x)-E\right)_+ + \left(E(t,x-\Delta_x)-E\right)_+\right]\right] \right\}\nonumber\\
&&+\Delta_x\frac{\partial n}{\partial x}\left\{ \left[ 1+\uptau R\right]\left[\frac{\theta}{2}+\frac{\eta \mu}{2 E_{\rm max}} \left(E-E(t,x+\Delta_x)\right)_+\right]\right\}-\Delta_x\frac{\partial n}{\partial x}\left\{ \left[ 1+\uptau R\right]\left[\frac{\theta}{2}+\frac{\eta \mu}{2 E_{\rm max}} \left(E-E(t,x-\Delta_x)\right)_+\right]\right\}\nonumber\\%
&&+\Delta_y\frac{\partial n}{\partial y}\left\{ \beta\left[ 1+\uptau R\right]\left[\frac{\theta}{2}+\frac{\eta \mu}{2 E_{\rm max}} \left(E-E(t,x+\Delta_x)\right)_+\right]\right\}-\Delta_y\frac{\partial n}{\partial y}\left\{ \beta\left[ 1+\uptau R\right]\left[\frac{\theta}{2}+\frac{\eta \mu}{2 E_{\rm max}} \left(E-E(t,x-\Delta_x)\right)_+\right]\right\}\nonumber\\%
&&+\frac{\Delta_x^2}{2}\frac{\partial^2 n}{\partial x^2}\left\{ \left[ 1+\uptau R\right]\left[\frac{\theta}{2}+\frac{\eta \mu}{2 E_{\rm max}} \left(E-E(t,x+\Delta_x)\right)_+\right]\right\}+\frac{\Delta_x^2}{2}\frac{\partial^2 n}{\partial x^2}\left\{ \left[ 1+\uptau R\right]\left[\frac{\theta}{2}+\frac{\eta \mu}{2 E_{\rm max}} \left(E-E(t,x-\Delta_x)\right)_+\right]\right\}\nonumber\\%
&&+\frac{\Delta_y^2}{2}\frac{\partial^2 n}{\partial y^2}\left\{ \beta\left[ 1+\uptau R\right]\left[\frac{\theta}{2}+\frac{\eta \mu}{2 E_{\rm max}} \left(E-E(t,x+\Delta_x)\right)_+\right]\right\}+\frac{\Delta_y^2}{2}\frac{\partial^2 n}{\partial y^2}\left\{ \beta\left[ 1+\uptau R\right]\left[\frac{\theta}{2}+\frac{\eta \mu}{2 E_{\rm max}} \left(E-E(t,x-\Delta_x)\right)_+\right]\right\}\nonumber\\%
&&+\frac{\Delta_y^2}{2}\frac{\partial^2 n}{\partial y^2}\left\{ \beta\left[ 1+\uptau R\right]\left[1-\theta-\frac{\eta \mu}{2 E_{\rm max}}\left[ \left(E(t,x+\Delta_x)-E\right)_+ + \left(E(t,x-\Delta_x)-E\right)_+\right]\right] \right\} +h.o.t.\nonumber
\end{eqnarray}

Cancelling out $\theta$ terms, where possible, we obtain

\begin{eqnarray}
\label{Der7}
&&n(t+\uptau,x,y)=n\left\{ \left[ 1+\uptau R\right]\left[\frac{\eta \mu}{2 E_{\rm max}} \left(E-E(t,x+\Delta_x)\right)_+\right]\right\}+n\left\{ \left[ 1+\uptau R\right]\left[\frac{\eta \mu}{2 E_{\rm max}} \left(E-E(t,x-\Delta_x)\right)_+\right]\right\}\\%
&&+n\left\{ \left[ 1+\uptau R\right]\left[1-\frac{\eta \mu}{2 E_{\rm max}}\left[ \left(E(t,x+\Delta_x)-E\right)_+ + \left(E(t,x-\Delta_x)-E\right)_+\right]\right] \right\}\nonumber\\
&&+\Delta_x\frac{\partial n}{\partial x}\left\{ \left[ 1+\uptau R\right]\left[\frac{\eta \mu}{2 E_{\rm max}} \left(E-E(t,x+\Delta_x)\right)_+\right]\right\}-\Delta_x\frac{\partial n}{\partial x}\left\{ \left[ 1+\uptau R\right]\left[\frac{\eta \mu}{2 E_{\rm max}} \left(E-E(t,x-\Delta_x)\right)_+\right]\right\}\nonumber\\%
&&+\Delta_y\frac{\partial n}{\partial y}\left\{ \beta\left[ 1+\uptau R\right]\left[\frac{\eta \mu}{2 E_{\rm max}} \left(E-E(t,x+\Delta_x)\right)_+\right]\right\}-\Delta_y\frac{\partial n}{\partial y}\left\{ \beta\left[ 1+\uptau R\right]\left[\frac{\eta \mu}{2 E_{\rm max}} \left(E-E(t,x-\Delta_x)\right)_+\right]\right\}\nonumber\\%
&&+\frac{\Delta_x^2}{2}\frac{\partial^2 n}{\partial x^2}\left\{ \left[ 1+\uptau R\right]\left[\theta+\frac{\eta \mu}{2 E_{\rm max}} \left(E-E(t,x+\Delta_x)\right)_+\right]\right\}+\frac{\Delta_x^2}{2}\frac{\partial^2 n}{\partial x^2}\left\{ \left[ 1+\uptau R\right]\left[\frac{\eta \mu}{2 E_{\rm max}} \left(E-E(t,x-\Delta_x)\right)_+\right]\right\}\nonumber\\%
&&+\frac{\Delta_y^2}{2}\frac{\partial^2 n}{\partial y^2}\left\{ \beta\left[ 1+\uptau R\right]\left[\frac{\eta \mu}{2 E_{\rm max}} \left(E-E(t,x+\Delta_x)\right)_+\right]\right\}+\frac{\Delta_y^2}{2}\frac{\partial^2 n}{\partial y^2}\left\{ \beta\left[ 1+\uptau R\right]\left[\frac{\eta \mu}{2 E_{\rm max}} \left(E-E(t,x-\Delta_x)\right)_+\right]\right\}\nonumber\\%
&&+\frac{\Delta_y^2}{2}\frac{\partial^2 n}{\partial y^2}\left\{ \beta\left[ 1+\uptau R\right]\left[1-\frac{\eta \mu}{2 E_{\rm max}}\left[ \left(E(t,x+\Delta_x)-E\right)_+ + \left(E(t,x-\Delta_x)-E\right)_+\right]\right] \right\} +h.o.t.\nonumber
\end{eqnarray}

We can then rearrange the equations to obtain
\begin{eqnarray}
\label{Der7b}
&&n(t+\uptau,x,y)=n\left\{ \left[ 1+\uptau R\right]\left[\frac{\eta \mu}{2 E_{\rm max}} \left(E-E(t,x+\Delta_x)\right)_+\right]\right\}+n\left\{ \left[ 1+\uptau R\right]\left[\frac{\eta \mu}{2 E_{\rm max}} \left(E-E(t,x-\Delta_x)\right)_+\right]\right\}\\%
&&+n\left[ 1+\uptau R\right]-n\left\{ \left[ 1+\uptau R\right]\left[\frac{\eta \mu}{2 E_{\rm max}}\left[ \left(E(t,x+\Delta_x)-E\right)_+\right]\right]\right\}-n\left\{ \left[ 1+\uptau R\right]\left[\frac{\eta \mu}{2 E_{\rm max}}\left[ \left(E(t,x-\Delta_x)-E\right)_+\right]\right]\right\}\nonumber\\
&&+\Delta_x\frac{\partial n}{\partial x}\left\{ \left[ 1+\uptau R\right]\left[\frac{\eta \mu}{2 E_{\rm max}} \left(E-E(t,x+\Delta_x)\right)_+\right]\right\}-\Delta_x\frac{\partial n}{\partial x}\left\{ \left[ 1+\uptau R\right]\left[\frac{\eta \mu}{2 E_{\rm max}} \left(E-E(t,x-\Delta_x)\right)_+\right]\right\}\nonumber\\%
&&+\Delta_y\frac{\partial n}{\partial y}\left\{ \beta\left[ 1+\uptau R\right]\left[\frac{\eta \mu}{2 E_{\rm max}} \left(E-E(t,x+\Delta_x)\right)_+\right]\right\}-\Delta_y\frac{\partial n}{\partial y}\left\{ \beta\left[ 1+\uptau R\right]\left[\frac{\eta \mu}{2 E_{\rm max}} \left(E-E(t,x-\Delta_x)\right)_+\right]\right\}\nonumber\\%
&&+\frac{\Delta_x^2}{2}\frac{\partial^2 n}{\partial x^2}\theta\left[ 1+\uptau R\right]+\frac{\Delta_x^2}{2}\frac{\partial^2 n}{\partial x^2}\left\{ \left[ 1+\uptau R\right]\left[\frac{\eta \mu}{2 E_{\rm max}} \left(E-E(t,x+\Delta_x)\right)_+\right]\right\}\nonumber\\
&&+\frac{\Delta_x^2}{2}\frac{\partial^2 n}{\partial x^2}\left\{ \left[ 1+\uptau R\right]\left[\frac{\eta \mu}{2 E_{\rm max}} \left(E-E(t,x-\Delta_x)\right)_+\right]\right\}\nonumber\\%
&&+\frac{\Delta_y^2}{2}\frac{\partial^2 n}{\partial y^2}\left\{ \beta\left[ 1+\uptau R\right]\left[\frac{\eta \mu}{2 E_{\rm max}} \left(E-E(t,x+\Delta_x)\right)_+\right]\right\}+\frac{\Delta_y^2}{2}\frac{\partial^2 n}{\partial y^2}\left\{ \beta\left[ 1+\uptau R\right]\left[\frac{\eta \mu}{2 E_{\rm max}} \left(E-E(t,x-\Delta_x)\right)_+\right]\right\}\nonumber\\%
&&+\frac{\Delta_y^2}{2}\frac{\partial^2 n}{\partial y^2}\beta\left[ 1+\uptau R\right]-\frac{\Delta_y^2}{2}\frac{\partial^2 n}{\partial y^2}\left\{ \beta\left[ 1+\uptau R\right]\left[\frac{\eta \mu}{2 E_{\rm max}}\left[ \left(E(t,x+\Delta_x)-E\right)_+\right]\right]\right\}\nonumber\\
&&-\frac{\Delta_y^2}{2}\frac{\partial^2 n}{\partial y^2}\left\{ \beta\left[ 1+\uptau R\right]\left[\frac{\eta \mu}{2 E_{\rm max}}\left[ \left(E(t,x-\Delta_x)-E\right)_+\right]\right]\right\} +h.o.t.\nonumber
\end{eqnarray}

Now, using the fact that, for real functions $f$, the relation $(f)_+-(-f)_+=f$ holds, we have
\begin{eqnarray}
\label{Der8}
&&n(t+\uptau,x,y)=n\left\{ \left[ 1+\uptau R\right]\left[\frac{\eta \mu}{2 E_{\rm max}} \left(E-E(t,x+\Delta_x)\right)\right]\right\}+n\left\{ \left[ 1+\uptau R\right]\left[\frac{\eta \mu}{2 E_{\rm max}} \left(E-E(t,x-\Delta_x)\right)\right]\right\}+n\left[ 1+\uptau R\right]\\
&&+\Delta_x\frac{\partial n}{\partial x}\left\{ \left[ 1+\uptau R\right]\left[\frac{\eta \mu}{2 E_{\rm max}} \left(E-E(t,x+\Delta_x)\right)_+\right]\right\}-\Delta_x\frac{\partial n}{\partial x}\left\{ \left[ 1+\uptau R\right]\left[\frac{\eta \mu}{2 E_{\rm max}} \left(E-E(t,x-\Delta_x)\right)_+\right]\right\}\nonumber\\%
&&+\Delta_y\frac{\partial n}{\partial y}\left\{ \beta\left[ 1+\uptau R\right]\left[\frac{\eta \mu}{2 E_{\rm max}} \left(E-E(t,x+\Delta_x)\right)_+\right]\right\}-\Delta_y\frac{\partial n}{\partial y}\left\{ \beta\left[ 1+\uptau R\right]\left[\frac{\eta \mu}{2 E_{\rm max}} \left(E-E(t,x-\Delta_x)\right)_+\right]\right\}\nonumber\\%
&&+\frac{\Delta_x^2}{2}\frac{\partial^2 n}{\partial x^2}\theta\left[ 1+\uptau R\right]+\frac{\Delta_x^2}{2}\frac{\partial^2 n}{\partial x^2}\left\{ \left[ 1+\uptau R\right]\left[\frac{\eta \mu}{2 E_{\rm max}} \left(E-E(t,x+\Delta_x)\right)_+\right]\right\}\nonumber\\
&&+\frac{\Delta_x^2}{2}\frac{\partial^2 n}{\partial x^2}\left\{ \left[ 1+\uptau R\right]\left[\frac{\eta \mu}{2 E_{\rm max}} \left(E-E(t,x-\Delta_x)\right)_+\right]\right\}\nonumber\\%
&&+\frac{\Delta_y^2}{2}\frac{\partial^2 n}{\partial y^2}\left\{ \beta\left[ 1+\uptau R\right]\left[\frac{\eta \mu}{2 E_{\rm max}} \left(E-E(t,x+\Delta_x)\right)\right]\right\}+\frac{\Delta_y^2}{2}\frac{\partial^2 n}{\partial y^2}\left\{ \beta\left[ 1+\uptau R\right]\left[\frac{\eta \mu}{2 E_{\rm max}} \left(E-E(t,x-\Delta_x)\right)\right]\right\}\nonumber\\
&&+\frac{\Delta_y^2}{2}\frac{\partial^2 n}{\partial y^2}\beta\left[ 1+\uptau R\right]+h.o.t.\nonumber
\end{eqnarray}
Next, assuming the function $E$ to be sufficiently regular, substituting the following Taylor expansion 
\[
E(t,x\pm\Delta_x)= E\pm \Delta_x\frac{\partial E}{\partial x}+\frac{\Delta_x^2}{2}\frac{\partial^2 E}{\partial x^2}+h.o.t.
\]
into \eqref{Der8} and removing higher order terms, we obtain
\begin{eqnarray}
\label{Der9}
&&n(t+\uptau,x,y)=n\left\{ \left[ 1+\uptau R\right]\left[\frac{\eta \mu}{2 E_{\rm max}} \left(-\Delta_x^2 \frac{\partial^2 E}{\partial x^2}\right)\right]\right\}+n\left[ 1+\uptau R\right]\\
&&+\Delta_x\frac{\partial n}{\partial x}\left\{ \left[ 1+\uptau R\right]\left[\frac{\eta \mu}{2 E_{\rm max}} \left(-\Delta_x \frac{\partial E}{\partial x}\right)_+\right]\right\}-\Delta_x\frac{\partial n}{\partial x}\left\{ \left[ 1+\uptau R\right]\left[\frac{\eta \mu}{2 E_{\rm max}} \left(\Delta_x \frac{\partial E}{\partial x}\right)_+\right]\right\}\nonumber\\%
&&+\Delta_y\frac{\partial n}{\partial y}\left\{ \beta\left[ 1+\uptau R\right]\left[\frac{\eta \mu}{2 E_{\rm max}} \left(-\Delta_x \frac{\partial E}{\partial x}\right)_+\right]\right\}-\Delta_y\frac{\partial n}{\partial y}\left\{ \beta\left[ 1+\uptau R\right]\left[\frac{\eta \mu}{2 E_{\rm max}} \left(\Delta_x \frac{\partial E}{\partial x}\right)_+\right]\right\}\nonumber\\%
&&+\frac{\Delta_x^2}{2}\frac{\partial^2 n}{\partial x^2}\theta\left[ 1+\uptau R\right]+\frac{\Delta_y^2}{2}\frac{\partial^2 n}{\partial y^2}\beta\left[ 1+\uptau R\right]+h.o.t.\nonumber
\end{eqnarray}
%
Absorbing terms of order $\mathcal{O}(\uptau \Delta_x^2)$ and $\mathcal{O}(\uptau \Delta_y^2)$ into \emph{h.o.t.} yields
\begin{eqnarray}
\label{Der10}
&&n(t+\uptau,x,y)=n\left\{ \left[\frac{\eta \mu}{2 E_{\rm max}} \left(-\Delta_x^2 \frac{\partial^2 E}{\partial x^2}\right)\right]\right\}+n\left[ 1+\uptau R\right]\\
&&+\Delta_x\frac{\partial n}{\partial x}\left\{\left[\frac{\eta \mu}{2 E_{\rm max}} \left(-\Delta_x \frac{\partial E}{\partial x}\right)_+\right]\right\}-\Delta_x\frac{\partial n}{\partial x}\left\{ \left[\frac{\eta \mu}{2 E_{\rm max}} \left(\Delta_x \frac{\partial E}{\partial x}\right)_+\right]\right\}\nonumber\\%
&&+\Delta_y\frac{\partial n}{\partial y}\left\{ \beta\left[\frac{\eta \mu}{2 E_{\rm max}} \left(-\Delta_x \frac{\partial E}{\partial x}\right)_+\right]\right\}-\Delta_y\frac{\partial n}{\partial y}\left\{ \beta\left[\frac{\eta \mu}{2 E_{\rm max}} \left(\Delta_x \frac{\partial E}{\partial x}\right)_+\right]\right\}\nonumber\\%
&&+\frac{\Delta_x^2}{2}\frac{\partial^2 n}{\partial x^2}\theta+\frac{\Delta_y^2}{2}\frac{\partial^2 n}{\partial y^2}\beta+h.o.t.\nonumber
\end{eqnarray}
Once again, using the relation $(f)_+-(-f)_+=f$ for real functions $f$, we find
\begin{eqnarray}
\label{Der11}
&&n(t+\uptau,x,y)=n\left\{ \left[\frac{\eta \mu}{2 E_{\rm max}} \left(-\Delta_x^2 \frac{\partial^2 E}{\partial x^2}\right)\right]\right\}+n\left[ 1+\uptau R\right]+\Delta_x\frac{\partial n}{\partial x}\left\{\left[\frac{\eta \mu}{2 E_{\rm max}} \left(-\Delta_x \frac{\partial E}{\partial x}\right)\right]\right\}\\
&&+\Delta_y\frac{\partial n}{\partial y}\left\{ \beta\left[\frac{\eta \mu}{2 E_{\rm max}} \left(-\Delta_x \frac{\partial E}{\partial x}\right)\right]\right\}+\frac{\Delta_x^2}{2}\frac{\partial^2 n}{\partial x^2}\theta+\frac{\Delta_y^2}{2}\frac{\partial^2 n}{\partial y^2}\beta+h.o.t.\nonumber
\end{eqnarray}
Further rearranging gives
\begin{eqnarray}
\label{Der12}
n(t+\uptau,x,y)=\frac{\theta\Delta_x^2 }{2}\frac{\partial^2 n}{\partial x^2}-\frac{\Delta_x^2 \eta \mu}{2 E_{\rm max}}\left[n\frac{\partial^2 E}{\partial x^2}+\frac{\partial n}{\partial x}\frac{\partial E}{\partial x}\right] +n+\uptau R n
+\frac{\Delta_y^2 \beta}{2}\frac{\partial^2 n}{\partial y^2}
- \beta\Delta_y \Delta_x\left[\frac{\eta \mu}{2 E_{\rm max}}\frac{\partial E}{\partial x}\right]\frac{\partial n}{\partial y}+h.o.t.
\end{eqnarray}
Dividing both sides of \eqref{Der12} by $\uptau$ yields
\begin{eqnarray}
&&\frac{n(t+\uptau,x,y)-n}{\uptau}=\frac{\theta\Delta_x^2 }{2\uptau}\frac{\partial^2 n}{\partial x^2}-\frac{\Delta_x^2 \eta \mu}{2 E_{\rm max}\uptau }\left[n\frac{\partial^2 E}{\partial x^2}+\frac{\partial n}{\partial x}\frac{\partial E}{\partial x}\right] + R n+\frac{\Delta_y^2 \beta}{2\uptau}\frac{\partial^2 n}{\partial y^2}- \frac{\beta\Delta_y \Delta_x}{\uptau}\left[\frac{\eta \mu}{2 E_{\rm max}}\frac{\partial E}{\partial x}\right]\frac{\partial n}{\partial y}+h.o.t.\nonumber
\end{eqnarray}
Now letting the time-step $\uptau\rightarrow 0$, the space-step $\Delta_x\rightarrow 0$ and the phenotype-step $\Delta_y\rightarrow0$ in such a way that
\[
\frac{\Delta_x^2 \theta}{2 \uptau}\rightarrow D \in \mathbb{R}^+_{*}, \quad \frac{\Delta_x^2 \eta}{2 E_{\rm max}\uptau}\rightarrow \nu \in \mathbb{R}^+_*,\quad \text{and} \quad\frac{\Delta_y^2 \beta}{2\uptau}\rightarrow \lambda \in \mathbb{R}^+_*
\]
and using the definition for $\chi(y)$ given by~\eqref{eq:defchimu}, we formally obtain 
\begin{eqnarray}
\label{Der14}
&&\frac{\partial n}{\partial t}=D\frac{\partial^2 n}{\partial x^2}-\chi(y)\left[n\frac{\partial^2 E}{\partial x^2}+\frac{\partial n}{\partial x}\frac{\partial E}{\partial x}\right] + R n+\lambda\frac{\partial^2 n}{\partial y^2},
\end{eqnarray}
which further simplifies to
\begin{eqnarray}
\label{Der15}
&&\frac{\partial n(t,x,y)}{\partial t}=D\frac{\partial^2 n(t,x,y)}{\partial x^2}-\chi(y)\frac{\partial}{\partial x}\left(n(t,x,y)\frac{\partial E(t,x)}{\partial x}\right)+R(y,\rho) n(t,x,y)+\lambda\frac{\partial^2 n(t,x,y)}{\partial y^2},\end{eqnarray}
from which, rearranging terms, we recover the PIDE~\eqref{eq:PDE}$_1$ for $n(t,x,y)$. 

{
\begin{remark}
Under the initial conditions and the assumptions on the model functions and in the asymptotic regime considered here, the cell density $\rho(t,x)$ behaves like a one-side compactly supported and monotonically decreasing travelling front, while the ECM density $E(t,x)$ is identically equal to $0$ on the support of $\rho(t,x)$ and identically equal to $1$ outside the support of $\rho(t,x)$ (i.e. ahead of the cell travelling front). Hence, since for each $t$ the relation ${\rm Supp}\left(n(t,x,y)\right) \subseteq {\rm Supp}\left(\rho(t,x)\right)$ holds $y$ by $y$, despite the fact that the ECM density $E(t,x)$ jumps from $0$ to $1$ moving from inside to outside  ${\rm Supp}\left(\rho(t,x)\right)$, we expect the formal method that we have employed to derive the PIDE~\eqref{Der15} from the underlying IB model to apply for $x \in {\rm Supp}\left(\rho(t,x)\right)$. This is confirmed by the good agreement between the cell dynamics predicted by numerical simulations of the IB model and numerical solutions of the corresponding continuum model. 
\end{remark}
}

\subsubsection{Equations for the MDE concentration and the ECM density}
Using first the fact that for $\uptau$, $\Delta_x$, and $\Delta_y$ sufficiently small the following relations hold
\begin{eqnarray}
n^{k}_{i,j}\approx n(t,x,y),\quad E^{k}_{i}\approx E(t,x),\quad M^{k}_{i}\approx M(t,x),\quad M^{k}_{i\pm1}\approx M(t,x\pm\Delta_x),\nonumber\\
n^{k+1}_{i,j}\approx n(t+\uptau,x,y),\quad E^{k+1}_{i}\approx E(t+\uptau,x),\quad M^{k+1}_{i}\approx M(t+\uptau,x),\quad p(y_j)\approx p(y),\nonumber
\end{eqnarray}
and then letting $\Delta_{\tau} \to 0$, $\Delta_{x} \to 0$, and $\Delta_y \to 0$ in the difference equations~\eqref{eq:MDE} and~\eqref{eq:ECM}, one formally obtains the PDE~\eqref{eq:PDE}$_3$ for $M(t,x)$ and the infinite-dimensional ODE~\eqref{eq:PDE}$_4$ for $E(t,x)$.
\end{document}